\documentstyle[11pt]{article}
       \font\tenmsb=msbm10
       \font\sevenmsb=msbm7
       \font\fivemsb=msbm5
       \catcode`\@=11
       \ifx\amstexloaded@\relax\catcode`\@=\active
       \endinput\else\let\amstexloaded@\relax\fi
       \def\spaces@{\space\space\space\space\space}
       \def\spaces@@{\spaces@\spaces@\spaces@\spaces@\spaces@}
       \def\space@.{\futurelet\space@\relax}
       \space@. %
       \def\Err@#1{\errhelp\defaulthelp@\errmessage{AmS-TeX error: #1}}
       \def\relaxnext@{\let\next\relax}
       \def\accentfam@{7}
       \def\tr{\triangle}
       \def\noaccents@{\def\accentfam@{0}}
       \def\Cal{\relaxnext@\ifmmode\let\next\Cal@\else
       \def\next{\Err@{Use \string\Cal\space only in math mode}}\fi\next}
       \def\Cal@#1{{\Cal@@{#1}}}
       \def\Cal@@#1{\noaccents@\fam\tw@#1}

       \def\Bbb{\relaxnext@\ifmmode\let\next\Bbb@\else
       \def\next{\Err@{Use \string\Bbb\space only in math mode}}\fi\next}
       
       \def\Bbb@#1{{\Bbb@@{#1}}}
       \def\Bbb@@#1{\noaccents@\fam\msbfam#1}
       \newfam\msbfam
       \textfont\msbfam=\tenmsb
       \scriptfont\msbfam=\sevenmsb
       \scriptscriptfont\msbfam=\fivemsb

\def\N{{\Bbb N}}
\def\Z{{\Bbb Z}}
\def\R{{\Bbb R}}
\def\T{{\Bbb T}}
\def\C{{\Bbb C}}

\newtheorem{Theorem}{Theorem}
\newtheorem{Lemma}{Lemma}[section]
\newtheorem{Proposition}{Proposition}

\@addtoreset{equation}{section}

\newcommand{\sss}{\smallskip}

\newcommand{\bs}{\bigskip}

\newcommand{\qed}{\nolinebreak\hfill\rule{2mm}{2mm}
\par\medbreak}
\newcommand{\proof}{\par\medbreak\it Proof: \rm}

\newcommand{\lep}{ \;<\;{\rm c}\; }
\newcommand{\la}{\langle }
\newcommand{\ra}{\rangle }

\newcommand{\kth}{e^{{\rm i}\langle k, \theta\rangle}\; }

\newcommand{\beq}{\begin{equation} }
\newcommand{\eeq}{\end{equation} }

\begin{document}
\setlength{\columnsep}{5pt}

\title{A  KAM Theorem for
Two-dimensional Nonlinear Schr\"odinger Equations\thanks{
The research was  supported by NNSFC Grant 11971012.}}
\author{\\  Jiansheng Geng, Shuaishuai Xue  \\
 {\footnotesize Department of Mathematics}\\
{\footnotesize Nanjing University, Nanjing 210093, P.R.China}\\
{\footnotesize Email: jgeng@nju.edu.cn; dg1621022@smail.nju.edu.cn}}


\date{}
\maketitle

\begin{abstract} We prove an infinite dimensional KAM theorem. As an application, we use the theorem to study
the two-dimensional nonlinear Schr\"{o}dinger
equation
$$iu_t-\triangle u +|u|^2u+\frac{\partial{f(x,u,\bar u)}}{\partial{\bar u}}=0, \quad t\in\Bbb R, x\in\Bbb T^2$$
with periodic boundary conditions, where the nonlinearity $\displaystyle f(x,u,\bar u)=\sum_{j,l,j+l\geq6}a_{jl}(x)u^j\bar u^l$, $a_{jl}=a_{lj}$ is a real analytic function in a neighborhood of the origin. We obtain for the equation  a Whitney smooth family of small--amplitude quasi--periodic
solutions which are partially hyperbolic.
\end{abstract}


\section{Introduction and Main Result}
The general problem discussed here is the persistency of quasi-periodic solutions of linear or integrable equations after Hamiltonian perturbation.
There have been many remarkable results in KAM
(Kolmogorov--Arnold--Moser) theory of Hamiltonian PDEs  achieved either by
methods from the  finite dimensional KAM theory
\cite{Ba3,
 CY,  E, EK, GXY, GY2, GY3, GY4, GY5,GY6, KaP, K1, K,
 KP,  P1, P2, P3, PP, PP2, W, XQY1, XQY2, Yuan2}, or by a Newtonian scheme
developed by  Craig, Wayne, Bourgain \cite{B1, B2, B3, B4, B5, B6,
BW, CW, Wa}. The advantage of the  method from the  finite dimensional
KAM theory is the construction of a local normal form in a
neighborhood of the obtained solutions in addition to the existence
of quasi-periodic solutions. The normal form is helpful to
understand the dynamics of the corresponding equations. For example, one sees the linear stability
and zero Lyapunov exponents. The scheme of  Craig-Wayne-Bourgain
avoids the cumbersome second Melnikov conditions by solving angle
dependent homological equations. The method is less Hamiltonian and more flexible than the KAM scheme to deal with resonant cases.
 All those  methods are well developed for one dimensional Hamiltonian PDEs.
 However, they meet difficulties in higher dimensional Hamiltonian PDEs.
 Bourgain \cite{B1} made the first breakthrough by proving that the two dimensional nonlinear Schr\"odinger equations
 admit small--amplitude quasi--periodic solutions.
 Later he improved in \cite{B4} his method and proved that the higher dimensional nonlinear Schr\"odinger and wave equations
 admit small--amplitude quasi--periodic solutions. Recently, W.--M.  Wang \cite{Wa} proved that the energy supercritical nonlinear Schr\"odinger equations admit small--amplitude quasi--periodic solutions.

 Constructing quasi-periodic solutions of higher dimensional
Hamiltonian PDEs  by  method developed from the  finite dimensional KAM theory appeared later. Geng--You \cite{GY3, GY4}
proved that  the higher dimensional nonlinear beam equations and
nonlocal Schr\"odinger  equations
 admit small--amplitude linearly--stable quasi--periodic solutions. The breakthrough of constructing quasi-periodic
solutions for
 more interesting higher dimensional Schr\"odinger equation by modified KAM method
 was made recently by
 Eliasson--Kuksin \cite{EK}. They proved that the higher dimensional nonlinear Schr\"odinger  equations
 admit small--amplitude linearly--stable quasi--periodic solutions. Quasi--periodic solutions of two dimensional cubic Schr\"odinger equation
    $$ {\rm i}u_t-\tr u+|u|^2 u=0,\qquad x\in \T^2,\ t\in \R, $$
  \noindent  with periodic boundary conditions are obtained by Geng--Xu--You \cite{GXY}. By carefully choosing tangential sites $\{i_1, \cdots,
i_b\}\in\Z^2$, the authors proved that the above nonlinear Schr\"odinger equation admits a family of
small-amplitude
 quasi-periodic solutions (see also \cite{PP2}). Very recently, Eliasson--Grebert--Kuksin \cite{EGK} proved that the higher dimensional nonlinear beam  equations
 admit small--amplitude  quasi--periodic solutions.

\sss   \sss   In this paper,our aim here is to pursue further  investigations  of the 2D nonlinear Schr\"odinger equation by developing the methods of \cite{GXY}. In \cite{GXY}, the authors require that the nonlinearity is independent of the space variable, so that a lot of technical complexity about the unbounded multiple  eigenvalues is successfully avoided. More precisely, when considering nonlinear Schr\"odinger equation especially in space dimension larger than one, a significant problem appears due to the presence of clusters of normal frequencies. Here the normal frequencies may have unbounded multiplicity
because the equation
\beq
m^2_1+m^2_2=R^2,m_1,m_2\in\Z\label{mmR}\eeq
(lattice points on a circle)
may have a large number of solutions for given $R$. It is important for our analysis that the integer solutions of (\ref{mmR}) appear in well-separated small clusters (of cardinality$\leq2$) and that the total number of integer solutions is at most $e^{\frac{\log R}{\log\log R}}\ll R^\varepsilon$.
The idea of the measure estimate  comes from
Geng--You\cite{GY6}. We use the
elementary repeated limit to substitute  Lipschitz
domain  by
 Eliasson--Kuksin \cite{EK},  thus our measure estimates are  easier and the whole proof is more KAM--like. More concretely, we consider the $2$-dimensional nonlinear
Schr\"{o}dinger equation
 \beq\label{nonlinearschro1} iu_t-\triangle u +|u|^2u+\frac{\partial{f(x,u,\bar u)}}{\partial{\bar u}}=0, \quad t\in\Bbb R, x\in\Bbb T^2
   \eeq
   with periodic boundary conditions
   $$
u(t,x_1+2\pi,x_2)=u(t,x_1,x_2+2\pi)=u(t,x_1,x_2),
$$
where $\displaystyle f(x,u,\bar u)=\sum_{j,l,j+l\geq6}a_{jl}(x)u^j\bar u^l,a_{jl}=a_{lj}$ is a real analytic function in a neighborhood of the origin.

 The operator
$A=-\triangle$ with periodic boundary conditions has
 eigenvalues $\{\lambda_n\}$ satisfying
$$\lambda_n=|n|^2=|n_1|^2+|n_2|^2, n=(n_1,n_2)\in \Z^2$$
 and the corresponding eigenfunctions
$\phi_n(x)=\frac{1}{2\pi}e^{{\rm i}\la n,x\ra}$ form a
basis in the domain of the operator.

A finite set
$S=\{i_1,\cdots,i_b\}\subset\Z^2$ is called \emph{admissible} if\\
\noindent
1.Any three of them are not vertices of a rectangle.\\
2.For any $n\in \Z^2\setminus S$, there exists at most one triplet $\{i,j,m\}$ with $i,j\in S$, $m\in \Z^2\setminus S$ such that $n-m+i-j=0$ and
$|n|^2-|m|^2+|i|^2-|j|^2=0$. If such triplet exists, we say that $n,m$ are resonant of the first type. By definition, $n,m$ are mutually uniquely determined. We say that $(n,m)$ is a resonant pair of first type. Geometrically, $(n,m,i,j)$ forms a rectangle with $n,m$ being two adjacent vertices.\\
3.For any $n\in \Z^2\setminus S$, there exists at most one triplet $\{i,j,m\}$ with $i,j\in S$, $m\in \Z^2\setminus S$ such that $n+m-i-j=0$ and
$|n|^2+|m|^2-|i|^2-|j|^2=0$. If such triplet exists, we say that $n,m$ are resonant of the second type. By definition, $n,m$ are mutually uniquely determined. We say that $(n,m)$ is a resonant pair of second type. Geometrically, $(n,m,i,j)$ forms a rectangle with $n,m$ being two diagonal vertices.\\
4.Any $n\in \Z^2\setminus S$ is not resonant of both the first type and the second type, i.e., there exist no $i,j,f,g\in S$ and $m,m'\in \Z^2\setminus S$, such that
$$ \left\{
\begin{array}{lcl}
n-m+i-j=0\\
|n|^2-|m|^2+|i|^2-|j|^2=0\\
n+m'-f-g=0\\
|n|^2+|m'|^2-|f|^2-|g|^2=0
\end{array} \right. $$
Geometrically, any two of the above defined rectangles cannot share vertex in $\Z^2\setminus S$.

In Appendix A of \cite{GXY}, a concrete way of constructing the admissible set is given. It is plausible that any randomly chosen set $S$ is almost surely admissible.

Now we state the main theorem as follows.
\begin{Theorem}\label{main}
Let $S=\{i_1,\cdots,i_b\}\subset\Z^2$ be an admissible set. There exists a Cantor set $\Cal C$ of positive--measure such that for any $\xi=(\xi_1,\cdots,\xi_b)\in\Cal C$, when $\xi_i^2+\xi_j^2<14\xi_i\xi_j$, the nonlinear Schr\"odinger equation (\ref{nonlinearschro1}) admits  a
small--amplitude,
 quasi--periodic solution of the form
\[ u(t, x)=\sum_{j=1}^b \sqrt{\xi_j}
e^{{\rm i}\omega_j t}\phi_{i_j}+O(|\xi|^{\frac32}),\omega_j=|i_j|^2+O(|\xi|).\]

\end{Theorem}

\sss \noindent {\bf Remark } We require $\xi_i^2+\xi_j^2<14\xi_i\xi_j$, such that the obtained tori are partially hyperbolic. When $\xi_i^2+\xi_j^2\geq 14\xi_i\xi_j$, one can prove the existence of the elliptic tori, however the proof is more complicated, which will be considered in the forthcoming
paper.

 This paper is organized as follows: In section 2 we give an infinite dimensional KAM theorem;
  in section 3, we give its application
 to two-dimensional Schr\"odinger equations.
  The proof of the KAM theorem is
  given in  section
 4, 5, 6. Some technical lemmas are given in the Appendix.

\section {An Infinite Dimensional KAM Theorem for Hamiltonian
Partial Differential Equations}

\sss In this section, we will formulate an infinite dimensional KAM
theorem that can be applied to two-dimensional  Schr\"odinger
equations under periodic boundary conditions.

We start by introducing some notations. For given $b$ vectors in
$\Z^2$, say $\{i_1, \cdots, i_b\}$, we denote
$\Z^2_1=\Z^2\setminus \{ i_1, \cdots, i_b\}$. Let $w=(\cdots,
w_n,\cdots)_{n\in \Z^2_1}$, and its complex conjugate $\bar
w=(\cdots,\bar w_n,\cdots)_{n\in \Z^2_1}$. We introduce the
weighted norm
$$\|w\|_{\rho} =\sum_{{n\in\Z^2_1}}|w_n|e^{|n|\rho},$$
where $|n|=\sqrt{n_1^2+n_2^2}$, $n=(n_1,n_2)\in \Z^2$ and
 $\rho
> 0$. Denote a neighborhood of $\T^b\times\{I=0\}\times\{w=0\}\times\{\bar w=0\}$ by
$$D_\rho(r,s)=\{(\theta,I,w,\bar w):|{\rm Im} \theta|<r,|I|<s^2,{\|w\|}_{\rho}<s,
{\|\bar w\|}_{\rho}<s\},$$
 where $|\cdot|$ denotes the sup-norm of complex vectors.
Moreover, we denote by $\Cal O$ a positive--measure parameter set
in $\R^b$.

 Let
$\alpha\equiv (\cdots,\alpha_n,\cdots)_{n\in\Z_1^2}$, $\beta
\equiv (\cdots, \beta _n, \cdots)_{n\in\Z_1^2}$,  $\alpha_n$ and $
\beta_n\in \N $ with finitely many non-zero  components of
positive integers. The product $w^{\alpha} \bar w^{\beta }$
denotes $\prod_n w_n^{\alpha_n}\bar w_n^{\beta_n}$. For any given
function \beq\label{2.2}
 F(\theta, I, w, \bar w)=\sum_{\alpha,\beta }
F_{\alpha\beta}(\theta, I)w^{\alpha} \bar w^{\beta },\eeq where
$\displaystyle F_{\alpha\beta}=\sum_{k\in\Z^b,l\in
\N^b}F_{kl\alpha\beta }(\xi)I^le^{{\rm i}
 \la k,\theta\ra}$ is $C_W^4$ function  in parameter $\xi$ in the
sense of Whitney, we denote \beq\label{2.4}\|F\|_{\mathcal O}=
\sum_{\alpha,\beta,k,l } |F_{kl\alpha \beta }|_{\mathcal O}\
|I^{l}|e^{|k||{\rm Im} \theta|}\,|w^{\alpha}||\bar w^{\beta }|\eeq where
 $|F_{kl\alpha \beta }|_{\mathcal O}$ is short for

 $$|F_{kl\alpha \beta }|_{\cal O}\equiv \sup_{\xi\in \Cal O}\sum_{0\leq d\leq 4}|{\partial_\xi^d F_{kl\alpha \beta }}|.$$ \noindent (the derivatives with respect to $\xi$ are in the sense of Whitney).
We define the weighted norm of $F$ by
 \begin{equation}\label{2.3}
  \|F\|_{ D_\rho(r,s)  ,\mathcal O}\equiv \sup_{D_\rho(r,s)}\|F\|_{\mathcal O},
\end{equation}

 To a function $F$, we associate a Hamiltonian vector
field defined by
$$ X_F=(F_I, -F_\theta, \{{\rm
i}F_{w_n}\}_{n\in \Z_1^2}, \{-{\rm i}F_{\bar
w_n}\}_{n\in\Z_1^2}).$$  Its
 weighted  norm is defined by \footnote{The norm  $\|\cdot\|_{D_\rho( r,s), \cal O}$ for scalar
functions is defined in (\ref{2.3}). The vector  function $G: D_\rho(
r,s)\times {\cal O}\to \C^m$, ($m<\infty$) is similarly defined as
$\|G\|_{D_\rho( r,s), \cal O}=\sum_{i=1}^m\|G_i\|_{D_\rho( r,s), \cal O}$.}
\begin{eqnarray}
\|X_F\|_{\!{}_{ D_\rho(r,s)  , \cal O}}&\equiv& \|F_I\|_{ D_\rho(r,s)  , \cal O}+ \frac
1{s^2}\|F_\theta\|_{ D_\rho(r,s)  , \cal O}\nonumber\\
&+&\sup_{D_\rho(r,s)}[ \frac 1s\sum_{n\in\Z_1^2} \|F_{w_n}\|_{\cal
O}e^{|n|\rho}+ \frac 1s\sum_{n\in\Z_1^2} \|F_{\bar
w_n}\|_{\cal O} e^{|n|\rho}
]\label{2.6}
\end{eqnarray}

Suppose that $S$ is an admissible set.
Let ${\cal L}_2$ be the subset of $Z^2_1$ with the following property: for each $n\in {\cal L}_2$,
there exists a unique triplet $(i,j,m)$ with $m\in Z^2_1$,$i,j\in S$ such that
$$-i-j+n+m=0,-|i|^2-|j|^2+|n|^2+|m|^2=0.$$

 We now describe a family of Hamiltonians studied in this paper. Let
 $$H_0=N+{\cal B}+\bar{\cal B},$$
\begin{eqnarray*}
N &=&\la\omega(\xi),I\ra+ \sum_{n\in \Z_1^2\backslash{{\cal L}_2}}\Omega_n(\xi)w_n \bar
w_n+\sum_{n'\in {{\cal L}_2}}(\Omega_{n'}(\xi)-\omega_{i'}(\xi))w_{n'}
\bar
w_{n'}
\end{eqnarray*}
Recall that $(i',j')$ is uniquely determined by the corresponding resonant pair $(n',m')$ in ${\cal L}_2$.
$${\cal B}=\sum_{n'\in {\cal L}_2}a_{n'}(\xi)w_{n'} w_{m'}$$
$$\bar{\cal B}=\sum_{n'\in {\cal L}_2}\bar a_{n'}(\xi)\bar w_{n'} \bar w_{m'}$$
where  $\xi\in \Cal O$ is a parameter, the phase space is endowed
with the symplectic structure $\displaystyle dI\wedge d\theta +
{\rm i} \sum_{n\in \Z_1^2} dw_n \wedge d \bar w_n$.

\sss For each $\xi\in \Cal O$, the Hamiltonian equation for $H_0$ admits
 special solutions $(\theta, 0, 0, 0)\to (\theta+\omega t,
0,0,0)$ that  corresponds to
 an invariant torus on the phase space.

\sss Consider now  the perturbed Hamiltonian
\begin{equation}\label{hamH}
H=H_0+P=N+{\cal B}+\bar{\cal B}+P(\theta,I,w,\bar w, \xi).
\end{equation}
Our goal is to prove that, for most values of parameter
 $\xi \in \Cal O$ (in Lebesgue measure
sense), the Hamiltonians $H=N+{\cal B}+\bar{\cal B}+P$ still admit
 invariant tori provided  that $\|X_P\|_{\!{}_{ D_\rho(r,s)  , \cal O}}$ is sufficiently
small.

\sss Decomposition of $\Z_1^2\backslash{{\cal L}_2}$.
For a nonnegative integer $\Delta$ we define an equivalence relation on $\Z_1^2\backslash{{\cal L}_2}$ generated by the pre-equivalence relation
$$a\sim b \Longleftrightarrow
\{
|a|^2=|b|^2 ,|a-b|\leq\Delta\} $$
Let $[a]_\Delta$ denote the equivalence class (block) and let $(\Z_1^2\backslash{{\cal L}_2})_\Delta$ be the set of equivalence classes. It is trivial that each block $[a]_\Delta$ is finite (we will write $[\cdot]$ for $[\cdot]_\bigtriangleup$).

Case 1:$|a|\leq \Delta$, we know $\sharp\{b:|a|=|b|,b\in \Z^2\}\leq e^{\frac{\log\Delta}{\log\log\Delta}}\ll \Delta^\varepsilon$;

Case 2:$|a|>\Delta$, we have $\sharp\{b:|a|=|b|,|a-b|\leq\Delta^{\frac{1}{3}},b\in \Z^2\}\leq2$.

  In order to have a compact formulation when solving homological equations, we rewrite $H$ into matrix form. Let $z_{[n]}=(w_i)_{i\in [n]}$, $\bar z_{[n]}=(\bar{w_i})_{i\in [n]}$; else $z_n=w_n,\bar{z}_n=\bar{w}_n$.
\begin{eqnarray*}
H&=&\la\omega(\xi),I\ra+ \sum_{n\in \Z_1^2\backslash{{\cal L}_2}}\Omega_n(\xi)w_n \bar
w_n+\sum_{n'\in {{\cal L}_2}}(\Omega_{n'}(\xi)-\omega_{i'}(\xi))w_{n'}
\bar
w_{n'}+{\cal B}+\bar{\cal B}+P\\
 &=&\la\omega(\xi),I\ra +\sum_{[n]} \la A_{[n]} z_{[n]}, \bar z_{[n]}\ra+\sum_{n'\in {{\cal L}_2}}(\Omega_{n'}(\xi)-\omega_{i'}(\xi))z_{n'}
\bar
z_{n'} +{\cal B}+\bar{\cal B}+P
\end{eqnarray*}
 where $A_{[n]}$ is $\sharp [n]\times\sharp [n]$  matrix.\\

We consider Hamiltonian $H$ satisfying the following hypotheses:

\bs \noindent $(A1)${\it Nondegeneracy:} The map $\xi\to
\omega(\xi)$ is a $C^4_W(\Cal O)$ diffeomorphism between $\Cal O$ and its
image.

\bs \noindent $(A2)${\it Asymptotics of normal frequencies:} \begin{equation}\label{asymp1} \Omega_n=\varepsilon^{-a}|n|^2+\tilde{\Omega}_n,a\geq0,n\in\Z_1^2\backslash{{\cal L}_2}
\end{equation}

where $\tilde{\Omega}_n$'s  are $C^4_W(\Cal O)$ functions of $\xi$ with $C^4_W(\Cal O)$-norm bounded by some  positive constant $L$.

\bs \noindent $(A3)$ {\it Melnikov's non--resonance conditions:} For $n\in\Z_1^2\backslash{{\cal L}_2}$, let
$$A_{[n]}=\Omega_{[n]}+(P_{i j }^{011})_{i\in [n],j\in [n]}=(\Omega_{i j }+P_{i j }^{011})_{i\in [n],j\in [n]}$$
where if $i\neq j$,$\Omega_{i j }=0$; if $i = j$, $\Omega_{i j }=\Omega_i$. When $|i-j|> K$, $P_{i j }^{011}=0$.
where ${A_{[n]} }'s$ are $C^4_W$ functions of $\xi$ with $C^4_W$-norm bounded by some  positive constant $L$, that is to say
$$\sup_{\xi\in \Cal O}\max_{0< d\leq4}\|{\partial^{d}_\xi A_{[n]}}\|\leq L$$
we assume that $\omega(\xi)$,$A_{[n]}(\xi)\in C^4_W(\Cal O)$ and there exist $\gamma,\tau>0$ such that, for $|k|\leq K$,

\[ |\langle k,\omega\rangle|\ge \frac{\gamma}{K^\tau}, k\neq 0,\]
 \[|\langle k,\omega\rangle \pm\widetilde{\lambda}_j|\ge \frac{\gamma}{K^\tau},j\in{[n]} \]
 \[|\langle k,\omega\rangle \pm \widetilde{\lambda}_i\pm \widetilde{\lambda}_j|\ge \frac{\gamma}{K^\tau},i\in{[m]},j\in{[n]} \]
where $\widetilde{\lambda}_i,\widetilde{\lambda}_j$ are $A_{[n]}$ and $A_{[m]}$'s eigenvalues respectively.

Let
$${\cal A}_n=A_{[n]},n\in\Z_1^2\backslash{{\cal L}_2}$$
$$ {\cal A}_n=\left(
  \begin{array}{cccc}
\Omega_{n}-\omega_{i} & -\frac{1}{2\pi^2}\sqrt{\xi_{i}\xi_{j}}\\
    \frac{1}{2\pi^2}\sqrt{\xi_{i}\xi_{j}}& -(\Omega_{m}-\omega_{j})
    \end{array}
\right) \\,n\in {\cal L}_2$$
where $(n,m)$ are resonant pairs, $(i,j)$
are uniquely determined by $(n,m)$ in ${\cal L}_2$.  \\
We assume that $\omega(\xi)$, ${\cal A}_n(\xi)\in C^4_W(\Cal O)$ and there exist $\gamma,\tau>0$ such that\footnote{The tensor product (or direct product) of two $m\times n,k\times l$ matrices $A=(a_{ij}),B$ is a $(mk)\times(nl)$ matrix defined by
$$ A\otimes B=(a_{ij}B)=\left(
  \begin{array}{cccc}
a_{11}B &\cdots& a_{1n}B\\
\cdots&\cdots&\cdots\\
    a_{m1}B&\cdots & a_{mn}B
    \end{array}
\right) \\\cdots$$
$\|\cdot\|$ for matrix denotes the operator norm, i.e., $\|M\|=\sup_{|y|=1}|My|$. Recall that $\omega$ and ${\cal A}_n,{\cal A}'_n$ depend on $\xi$.}
 (here $I_2$ is $2\times2$ identity matrix)
$$|det(\la k,\omega\ra I\pm{\cal A}_n\otimes I_2 \pm I_2\otimes {\cal A}_{n'})|\geq \frac{\gamma}{K^\tau},
k\neq0,n,n'\in{\cal L}_2.$$
{\bf We assume that the eigenvalues of ${\cal A}_n$ ($n\in {\cal L}_2$) have the non--zero imaginary parts so that the obtained tori are partially hyperbolic.}

 \bs \noindent $(A4)$ {\it Regularity of ${\cal B}+\bar{\cal B}+P$:} ${\cal B}+\bar{\cal B}+P$ is real analytic in $I,\theta,w,\bar w$ and Whitney smooth in $\xi$; in addition
$$\|X_{\cal B}\|_{D_\rho(r,s), \cal O}<1,\|X_P\|_{D_\rho(r,s), \cal O}<\varepsilon$$

\bs \noindent $(A5)$ {\it T\"{o}plitz-Lipschitz property:} For any
fixed $n, m\in {\Z}^2$, $c\in\Z^2\setminus\{0\}$, the limits
$$\lim_{t\to \infty}\frac{\partial^2 ({\cal B}+P)}{\partial
w_{n+tc}\partial w_{m-tc}}, \quad \lim_{t\to \infty}\frac{\partial^2 (\sum_{n\in {\Z}_1^2}\tilde{\Omega}_nw_n\bar w_n+P)}{\partial w_{n+tc}\partial \bar w_{m+tc}},\quad
\lim_{t\to \infty}\frac{\partial^2 (\bar{\cal B}+P)}{\partial \bar w_{n+tc}\partial \bar w_{m-tc}}$$ exist. Moreover,
 there exists $K>0$, such that when
$|t|>K$, $N+{\cal B}+\bar{\cal B}+P$ satisfies
 \[ \|\frac{\partial^2 ({\cal B}+P)}{\partial w_{n+tc}\partial w_{m-tc}}
 -\lim_{t\to \infty}\frac{\partial^2 ({\cal B}+P)}{\partial w_{n+tc}\partial w_{m-tc}}
 \|_{\!{}_{D_\rho(r,s), \Cal
O}}\le \frac{\varepsilon}{|t|}e^{-|n+m|\rho},
 \]
 \[ \|\frac{\partial^2 (\sum_{n\in {\Z}_1^2}\tilde{\Omega}_nw_n\bar w_n+P)}{\partial w_{n+tc}\partial \bar w_{m+tc}}
 -\lim_{t\to \infty}
 \frac{\partial^2 (\sum_{n\in {\Z}_1^2}\tilde{\Omega}_nw_n\bar w_n+P)}
 {\partial w_{n+tc}\partial \bar w_{m+tc}}\|_{\!{}_{D_\rho(r,s), \Cal
O}}\le \frac{\varepsilon}{|t|}e^{-|n-m|\rho},
 \]
 \[ \|\frac{\partial^2 (\bar{\cal B}+P)}{\partial \bar w_{n+tc}\partial \bar w_{m-tc}}-
 \lim_{t\to \infty}\frac{\partial^2 (\bar{\cal B}+P)}{\partial \bar w_{n+tc}\partial \bar w_{m-tc}}
 \|_{\!{}_{D_\rho(r,s), \Cal
O}}\le \frac{\varepsilon}{|t|}e^{-|n+m|\rho}.
 \]

\bs \noindent Now we are ready to state an infinite dimensional KAM Theorem.

\begin{Theorem}\label{KAM}
Assume that the Hamiltonian $H_0+P$ in (\ref{hamH}) satisfies
$(A1)$--$(A5)$. Let $\gamma>0$ be small enough, there exists a
positive constant $\varepsilon=\varepsilon(b, K,
\tau,\gamma,r,s,\rho)$.  Such that if $\|X_P\|_{\!{}_{D_\rho(r,s),
\cal O}}<\varepsilon$, then the following holds true: There exist a
Cantor set $\Cal O_\gamma\subset\Cal O$ with ${\rm meas}(\Cal
O\setminus \Cal O_\gamma)=O(\gamma^{\frac14})$ and two maps ( analytic in
$\theta$ and $C_W^4$ in $\xi$)
$$\Psi: \T^b\times \Cal O_\gamma\to D_\rho(r,s),\ \ \ \
 \tilde\omega:\Cal O_\gamma\to \R^b,$$ where $\Psi$ is
 $\frac{\varepsilon}{\gamma^4}$-close to the trivial embedding
$\Psi_0:\T^b\times \Cal O\to \T^b\times\{0,0,0\}$ and $\tilde
\omega$ is $\varepsilon$-close to the unperturbed frequency
$\omega$. Then for any $\xi\in \Cal O_\gamma$ and $\theta\in \T^b$,
the curve $t\to \Psi(\theta+\tilde\omega(\xi) t,\xi)$ is a
quasi-periodic solution of the Hamiltonian equations governed by
$H=H_0+P$. The obtained tori are partially hyperbolic.
\end{Theorem}

\section{Application to the Two-dimensional Schr\"odinger Equations}

\sss \noindent  We
 consider the two--dimensional nonlinear
Schr\"odinger equations
   \beq\label{nonlinearschro} {\rm i}u_t-\Delta u+|u|^2u+\frac{\partial{f(x,u,\bar u)}}{\partial{\bar u}}=0,
   \qquad x\in \T^2,\ t\in \R
   \eeq
   with periodic boundary conditions
   $$
u(t,x_1+2\pi,x_2)=u(t,x_1,x_2+2\pi)=u(t,x_1,x_2),
$$
where $\displaystyle f(x,u,\bar u)=\sum_{j,l,j+l\geq6}a_{jl}(x)u^j\bar u^l,a_{jl}=a_{lj}$ is a real analytic function in a neighborhood of the origin.

The operator
$A=-\triangle$ with periodic boundary conditions has
 eigenvalues $\{\lambda_n\}$ satisfying
$$\lambda_n=|n|^2=|n_1|^2+|n_2|^2, n=(n_1,n_2)\in \Z^2$$
 and the corresponding eigenfunctions
$\phi_n(x)=\frac{1}{2\pi}e^{{\rm i}\la n,x\ra}$ form a
basis in the domain of the operator.

  Equation (\ref{nonlinearschro}) can be
rewritten
 as a Hamiltonian equation
\beq\label{3.7+7} u_t={\rm i}\frac{\partial H}{\partial \bar u}
\eeq and the corresponding Hamiltonian is \beq\label{hamiltoniann}
H=\la Au,u\ra +\frac{1}{2}\int_{\T^2} |u|^4\ dx+\int_{\T^2} f(x,u,\bar u)\ dx,\eeq where
$\la \cdot,\cdot\ra$ denotes the inner product in $L^2$.

 Let
\[u(x)=\sum_{n\in \Z^2}{q_n}\phi_n(x),
\]
 System (\ref{3.7+7}) is then equivalent to the lattice
Hamiltonian equations
\begin{equation}\label{3.8+8} \dot q_n={\rm i}(\lambda_n q_n+  \frac{\partial G}{\partial \bar q_n}),
 \quad G\equiv \frac{1}{8\pi^2}\sum_{i-j+n-m=0}q_i\bar q_jq_n\bar q_m+\int_{\T^2} f(x,u,\bar u)\ dx ,\end{equation} with corresponding Hamiltonian
function
\begin{eqnarray}
H&=& \sum_{n\in \Z^2}\lambda_nq_n\bar q_n+\frac{1}{8\pi^2}\sum_{i-j+n-m=0}q_i\bar q_jq_n\bar q_m+\int_{\T^2} f(x,\sum_{n\in \Z^2}{q_n}\phi_n(x),\sum_{n\in \Z^2}{\bar{q}_n}\bar{\phi}_n(x))\ dx\nonumber\\
 &= &\sum_{n\in \Z^2}\lambda_n|q_n|^2+G\label{PH}\\
G&=& \frac{1}{8\pi^2}\sum_{i-j+n-m=0}q_i\bar q_jq_n\bar q_m+\int_{\T^2} f(x,\sum_{n\in \Z^2}{q_n}\phi_n(x),\sum_{n\in \Z^2}{\bar{q}_n}\bar{\phi}_n(x))\ dx\nonumber
\end{eqnarray}

 As in \cite {KP,P1,GY2}, the perturbation $G$ in
 (\ref{3.8+8})
has  the following  regularity property.
\begin{Lemma}\label{regularityGG}
For any
 fixed  $\rho>0$,
  the gradient $G_{\bar q}$ is real analytic as a map
 in a neighborhood of the origin with
 \beq \|G_{\bar q}\|_{\rho}\le c\|q\|_{\rho}^3.
\label{3.16+6} \eeq
\end{Lemma}
\proof
\begin{eqnarray*}
\|G_{\bar q}\|_{\rho}&=&\sum_{n\in\Z^2}|G_{\bar q_n}|e^{|n|\rho}\\
&\leq&c\sum_{{n,\alpha,\beta-e_n,|\alpha|+|\beta-e_n|=
3}}|q^\alpha\bar
q^{\beta-e_n}|e^{|n|\rho}\\
&\leq&c\sum_{\alpha,\beta-e_n,|\alpha|+|\beta-e_n|=
3}|q^\alpha\bar
q^{\beta-e_n}|e^{|\alpha|\rho}e^{|\beta-e_n|\rho}\\
&\leq&c\|q\|_{\rho}^3.
\end{eqnarray*}\qed

For an admissible set of tangential site
$S=\{i_1,\cdots,i_b\}\subset\Z^2$, we have a nice normal form for $H$.\\

\begin{Proposition}\label{P1}Let $S$ be admissible.For Hamiltonian function (\ref{PH}), there is a symplectic transformation $\Psi$, such that
\begin{equation}\label{P11}H\circ \Psi=\la \omega,I\ra+\la \Omega w,w\ra+{\cal A}+{\cal B}+\bar{\cal B}+P
\end{equation}
with
$$ \left\{
\begin{array}{lcl}
\omega_i(\xi)=\displaystyle\varepsilon^{-3}|i|^2-\frac{1}{4\pi^2}\xi_i+\sum_{j\in S}\frac{1}{2\pi^2}\xi_j\\
\Omega_n=\displaystyle\varepsilon^{-3}|n|^2+\sum_{j\in S}\frac{1}{2\pi^2}\xi_j
\end{array} \right. $$
$${\cal A}=\frac{1}{2\pi^2}\sum_{n\in {\cal L}_1}\sqrt{\xi_i\xi_j}w_n\bar w_me^{i\theta_i-i\theta_j}$$
$${\cal B}=\frac{1}{2\pi^2}\sum_{n'\in {\cal L}_2}\sqrt{\xi_{i'}\xi_{j'}}w_{n'} w_{m'}e^{-i\theta_{i'}-i\theta_{j'}}$$
$$\bar{\cal B}=\frac{1}{2\pi^2}\sum_{n'\in {\cal L}_2} \sqrt{\xi_{i'}\xi_{j'}}\bar w_{n'} \bar w_{m'}e^{i\theta_{i'}+i\theta_{j'}}$$
\begin{eqnarray}
|P|=&&O(\varepsilon^2|I|^2+\varepsilon^2|I|\|w\|^2_\rho+\varepsilon\xi^{\frac{1}{2}}\|w\|^3_\rho+\varepsilon^2\|w\|^4_\rho+\varepsilon\xi^3\nonumber\\
&&+\varepsilon^2\xi^{\frac{5}{2}}\|w\|_\rho+\varepsilon^3\xi^2\|w\|^2_\rho+\varepsilon^4\xi^{\frac{3}{2}}\|w\|^3_\rho).
\end{eqnarray}

\end{Proposition}

\proof
The proof consists of several sympiectic change of variables. Firstly, let
\begin{equation}\label{P12}
F=\sum_{{{i-j+n-m=0}\atop{|i|^2-|j|^2+|n|^2-|m|^2\neq0}}\atop\sharp S\cap\{i,j,n,m\}\geq2}\frac{i}{8\pi^2(\lambda_i-\lambda_j+\lambda_n-\lambda_m)}q_i\bar q_jq_n\bar q_m,
\end{equation}
and ${X}^1_F$ be the time one map of the flow of the associated Hamiltonian systems. The change of variables ${X}^1_F$ sends $H$ to

\begin{eqnarray}
H\circ{X}^1_F&=&H+ \{H,F\}+\int_0^1 (1-t)\{\{H,F\},F\}\circ \phi_F^{t}dt\nonumber\\
&=&\sum_{i\in S}\lambda_i|q_i|^2+\sum_{i\in \Z^2_1}\lambda_i|w_i|^2+\sum_{i\in S}\frac{1}{8\pi^2}|q_i|^4\\
&+&\sum_{i,j\in S,i\neq j}\frac{1}{2\pi^2}|q_i|^2|q_j|^2+\sum_{i\in S,j\in \Z^2_1}\frac{1}{2\pi^2}|q_i|^2|w_j|^2\\
&+&\sum_{n\in {\cal L}_1}\frac{1}{2\pi^2}q_i\bar q_jw_n\bar w_m+\sum_{{n'}\in {\cal L}_2}\frac{1}{2\pi^2}(q_{i'}q_{j'}\bar w_{n'}\bar w_{m'}+\bar q_{i'}\bar q_{j'} w_{n'} w_{m'})\label{P13}\\
&+&O(|q|\|w\|^3_\rho+\|w\|^4_\rho+|q|^6+|q|^5\|w\|_\rho+|q|^4\|w\|_\rho^2+|q|^3\|w\|_\rho^3).\nonumber
\end{eqnarray}
We remind that $(n,m)$ are resonant pairs and $(i,j)$ is uniquely determined by $(n,m)$; $(n',m')$ are resonant pairs and $(i',j')$ is uniquely determined by $(n',m')$ in (\ref{P13}).

  Next we introduce standard  action-angle variables in the tangential space
  $$q_j=\sqrt{I_j+\xi_j}e^{{\rm i}\theta_j},\bar q_j=\sqrt{I_j+\xi_j}e^{-{\rm i}\theta_j},j\in S,$$
  and
  $$q_n=w_n, \bar q_n=\bar w_n, n\in \Z^2_1,$$ we have
\begin{eqnarray*}
H\circ{X}^1_F=&&\sum_{i\in S}\lambda_i(I_i+\xi_i)+\sum_{i\in \Z^2_1}\lambda_i|w_i|^2+\sum_{i\in S}\frac{1}{8\pi^2}(I_i+\xi_i)^2\\
&+&\frac{1}{2\pi^2}\sum_{i,j\in S,i\neq j}(I_i+\xi_i)(I_j+\xi_j)+\frac{1}{2\pi^2}\sum_{i\in S,j\in \Z^2_1}(I_i+\xi_i)|w_j|^2\\
&+&\frac{1}{2\pi^2}\sum_{n\in {\cal L}_1}\sqrt{(I_i+\xi_i)(I_j+\xi_j)}w_n\bar w_me^{i\theta_i-i\theta_j}\\
&+&\frac{1}{2\pi^2}\sum_{{n'}\in {\cal L}_2}\sqrt{(I_{i'}+\xi_{i'})(I_{j'}+\xi_{j'})}w_{n'} w_{m'}e^{-i\theta_{i'}-i\theta_{j'}}\\
&+&\frac{1}{2\pi^2}\sum_{{n'}\in {\cal L}_2}\sqrt{(I_{i'}+\xi_{i'})(I_{j'}+\xi_{j'})}\bar w_{n'}\bar w_{m'}e^{i\theta_{i'}+i\theta_{j'}}\\
&+&O(\xi^{\frac12}\|w\|^3_\rho+\|w\|^4_\rho+\xi^3+\xi^{\frac52}\|w\|_\rho+\xi^2\|w\|_\rho^2
+\xi^{\frac32}\|w\|_\rho^3).\\
=&&\sum_{i\in S}\lambda_iI_i+\sum_{i\in \Z^2_1}\lambda_i|w_i|^2+\sum_{i\in S}\frac{1}{4\pi^2}\xi_i I_i
+\sum_{i,j\in S,i\neq j}\frac{1}{2\pi^2}\xi_iI_j+\sum_{i\in S,j\in \Z^2_1}\frac{1}{2\pi^2}\xi_i|w_j|^2\\
&+&\frac{1}{2\pi^2}\sum_{n\in {\cal L}_1}\sqrt{\xi_i\xi_j}w_n\bar w_me^{i\theta_i-i\theta_j}\\
&+&\frac{1}{2\pi^2}\sum_{n'\in {\cal L}_2}\sqrt{\xi_{i'}\xi_{j'}}w_{n'} w_{m'}e^{-i\theta_{i'}-i\theta_{j'}}\\
&+&\frac{1}{2\pi^2}\sum_{n'\in {\cal L}_2} \sqrt{\xi_{i'}\xi_{j'}}\bar w_{n'} \bar w_{m'}e^{i\theta_{i'}+i\theta_{j'}}\\
&+&O(|I|^2+|I|\|w\|^2_\rho+\xi^{\frac{1}{2}}\|w\|^3_\rho+\|w\|^4_\rho+\xi^3+\xi^{\frac{5}{2}}\|w\|_\rho+\xi^2\|w\|^2_\rho+\xi^{\frac{3}{2}}\|w\|^3_\rho)\\
=&&N+{\cal A}+{\cal B}+\bar{\cal B}+P
\end{eqnarray*}
where
$$N=\sum_{i\in S}\lambda_iI_i+\sum_{j\in \Z^2_1}\lambda_j|w_j|^2-\sum_{i\in S}\frac{1}{4\pi^2}\xi_i I_i
+\sum_{i,j\in S}\frac{1}{2\pi^2}\xi_iI_j+\sum_{i\in S,j\in \Z^2_1}\frac{1}{2\pi^2}\xi_i|w_j|^2$$
$${\cal A}=\frac{1}{2\pi^2}\sum_{n\in {\cal L}_1}\sqrt{\xi_i\xi_j}w_n\bar w_me^{i\theta_i-i\theta_j}$$
$${\cal B}=\frac{1}{2\pi^2}\sum_{n'\in {\cal L}_2}\sqrt{\xi_{i'}\xi_{j'}}w_{n'} w_{m'}e^{-i\theta_{i'}-i\theta_{j'}}$$
$$\bar{\cal B}=\frac{1}{2\pi^2}\sum_{n'\in {\cal L}_2} \sqrt{\xi_{i'}\xi_{j'}}\bar w_{n'} \bar w_{m'}e^{i\theta_{i'}+i\theta_{j'}}$$

By the scaling in time
$$\xi\rightarrow \varepsilon^3\xi,I\rightarrow \varepsilon^5I,\theta\rightarrow \theta,w\rightarrow\varepsilon^{\frac52}w,\bar w\rightarrow\varepsilon^{\frac52}\bar w$$
we finally arrive at the rescaled Hamiltonian
$$H=\varepsilon^{-8}H(\varepsilon^3\xi,\varepsilon^5I,\theta,\varepsilon^{\frac52}w,\varepsilon^{\frac52}\bar w)=\la \omega,I\ra+\la \Omega w,w\ra+{\cal A}+{\cal B}+\bar{\cal B}+P$$
where
$$ \left\{
\begin{array}{lcl}
\omega_i(\xi)=\displaystyle\varepsilon^{-3}|i|^2-\frac{1}{4\pi^2}\xi_i+\sum_{j\in S}\frac{1}{2\pi^2}\xi_j\\
\Omega_n=\displaystyle\varepsilon^{-3}|n|^2+\sum_{j\in S}\frac{1}{2\pi^2}\xi_j
\end{array} \right. $$
$${\cal A}=\frac{1}{2\pi^2}\sum_{n\in {\cal L}_1}\sqrt{\xi_i\xi_j}w_n\bar w_me^{i\theta_i-i\theta_j}$$
$${\cal B}=\frac{1}{2\pi^2}\sum_{n'\in {\cal L}_2}\sqrt{\xi_{i'}\xi_{j'}}w_{n'} w_{m'}e^{-i\theta_{i'}-i\theta_{j'}}$$
$$\bar{\cal B}=\frac{1}{2\pi^2}\sum_{n'\in {\cal L}_2} \sqrt{\xi_{i'}\xi_{j'}}\bar w_{n'} \bar w_{m'}e^{i\theta_{i'}+i\theta_{j'}}$$
\begin{eqnarray}
|P|=&&O(\varepsilon^2|I|^2+\varepsilon^2|I|\|w\|^2_\rho+\varepsilon\xi^{\frac{1}{2}}\|w\|^3_\rho+\varepsilon^2\|w\|^4_\rho+\varepsilon\xi^3\nonumber\\
&&+\varepsilon^2\xi^{\frac{5}{2}}\|w\|_\rho+\varepsilon^3\xi^2\|w\|^2_\rho+\varepsilon^4\xi^{\frac{3}{2}}\|w\|^3_\rho).
\end{eqnarray}\qed

We will show that, by a nonlinear symplectic coordinates transformation, the normal form in Proposition \ref{P1} can be transformed into the more elegant form. For this purpose, we need the following lemma from \cite{XY}.
\begin{Lemma}\label{symchange}
For any $k_1,k_2,\cdots,k_m\in\Z^b$, non-singular $m\times m$ matrix $S$ with $S^T\bar S=I$, the map $\Phi_0:(\theta,I,w,w)\rightarrow(\theta_+,I_+,z,z)$ defined by
$$ \left\{
\begin{array}{lcl}
\theta_+=\theta\\
I_+=I-\sum_{j=1}^{m}w_j\bar w_jk_j\\
z=SEw\\
\bar z=\bar S\bar E\bar w
\end{array} \right. $$
is symplectic with diagonal matrix
$$E=E(k_1,k_2,\cdots,k_m)=diag(e^{i\la k_1,\theta\ra},e^{i\la k_2,\theta\ra},\cdots,e^{i\la k_m,\theta\ra}).$$
\end{Lemma}
 The proof of the above lemma  refers to \cite{XY}.\\

A nonlinear symplectic coordinates transformation $\Phi$($\exists S$):
$$ \left\{
\begin{array}{lcl}
\theta_+=\theta\\
I_+=I-\displaystyle\sum_{n\in {\cal L}_1}(w_n\bar w_ne_i+w_m\bar w_me_j)+\sum_{n'\in {\cal L}_2}(w_{n'}\bar w_{n'}e_{i'}+w_{m'}\bar w_{m'}e_{j'})\\
{\left(
  \begin{array}{cccc}
z_n\\
z_m
  \end{array}
\right)}=S{\left(
  \begin{array}{cccc}
e^{i\la k_i,\theta\ra} & 0 \\
    0& e^{i\la k_j,\theta\ra}
  \end{array}
\right)}{\left(
  \begin{array}{cccc}
w_n\\
w_m
  \end{array}
  \right)},
  {\left(
  \begin{array}{cccc}
\bar z_n\\
\bar z_m
  \end{array}
\right)}=\bar S{\left(
  \begin{array}{cccc}
e^{-i\la k_i,\theta\ra} & 0 \\
    0& e^{-i\la k_j,\theta\ra}
  \end{array}
\right)}{\left(
  \begin{array}{cccc}
\bar w_n\\
\bar w_m
  \end{array}
  \right)},n\in {\cal L}_1\\
z_{n'}=w_{n'}e^{-i\theta_{i'}},\bar z_{n'}=\bar w_{n'}e^{i\theta_{i'}};z_{m'}=w_{m'}e^{-i\theta_{j'}},\bar z_{m'}=\bar w_{m'}e^{i\theta_{j'}},n'\in {\cal L}_2\\
z_n=w_n,\bar z_n=\bar w_n,n\in \Z^2_1\setminus({\cal L}_1\cup{\cal L}_2)
\end{array} \right. $$

We get Hamiltonian systems with the Hamiltonian
\begin{eqnarray}\label{Hamiltoniann}
H\circ\Psi\circ\Phi&=&\la\omega(\xi), I_+\ra+\sum_{{n\in \Z^2_1\setminus({\cal L}_1\cup{\cal L}_2)}} \Omega_n(\xi)z_n\bar
z_n\nonumber\\
&+&\sum_{{n\in {\cal L}_1}}[ (\varepsilon^{-3}(|n|^2+|i|^2)+\sum_{j\in S}\frac{1}{\pi^2}\xi_j-\frac{1}{8\pi^2}(\xi_i+\xi_j)+\frac{1}{8\pi^2}\sqrt{\xi^2_i+14\xi_i\xi_j+\xi_j^2})z_n\bar
z_n\nonumber\\
&+&(\varepsilon^{-3}(|m|^2+|j|^2)+\sum_{j\in S}\frac{1}{\pi^2}\xi_j-\frac{1}{8\pi^2}(\xi_i+\xi_j)-\frac{1}{8\pi^2}\sqrt{\xi^2_i+14\xi_i\xi_j+\xi_j^2})z_m\bar
z_m]\nonumber\\
&+&\sum_{{n'\in {\cal L}_2}}[ (\Omega_{n'}-\omega_{i'})z_{n'}\bar
z_{n'}+(\Omega_{m'}-\omega_{j'})z_{m'}\bar
z_{m'}]\nonumber\\
&+&\frac{1}{2\pi^2}\sum_{n'\in {\cal L}_2}\sqrt{\xi_{i'}\xi_{j'}}z_{n'} z_{m'}+\frac{1}{2\pi^2}\sum_{n'\in {\cal L}_2}\sqrt{\xi_{i'}\xi_{j'}}\bar z_{n'} \bar z_{m'}\nonumber\\
 &+&P(\theta_+, I_+, z, \bar z, \xi)\nonumber\\
 =&&N+{\cal B}+\bar{\cal B}+P
\end{eqnarray}

where

\begin{eqnarray*}N&=&\la\omega(\xi), I_+\ra+\sum_{{n\in \Z^2_1\setminus{\cal L}_2}} \Omega_n(\xi)z_n\bar
z_n\\
&+&\sum_{{n'\in {\cal L}_2}}[ (\Omega_{n'}-\omega_{i'})z_{n'}\bar
z_{n'}+(\Omega_{m'}-\omega_{j'})z_{m'}\bar
z_{m'}]
\end{eqnarray*}
$$ \left\{
\begin{array}{lcl}
\omega_i(\xi)=\displaystyle\varepsilon^{-3}|i|^2-\frac{1}{4\pi^2}\xi_i+\sum_{j\in S}\frac{1}{2\pi^2}\xi_j\\
\Omega_n=\displaystyle\varepsilon^{-3}|n|^2+\sum_{j\in S}\frac{1}{2\pi^2}\xi_j,{n\in \Z^2_1\setminus{\cal L}_1}\\
\Omega_n=\displaystyle\varepsilon^{-3}(|n|^2+|i|^2)+\sum_{j\in S}\frac{1}{\pi^2}\xi_j-\frac{1}{8\pi^2}(\xi_i+\xi_j)+\frac{1}{8\pi^2}\sqrt{\xi^2_i+14\xi_i\xi_j+\xi_j^2},n\in{\cal L}_1\\
\Omega_m=\varepsilon^{-3}(|m|^2+|j|^2)+\sum_{j\in S}\frac{1}{\pi^2}\xi_j-\frac{1}{8\pi^2}(\xi_i+\xi_j)-\frac{1}{8\pi^2}\sqrt{\xi^2_i+14\xi_i\xi_j+\xi_j^2},n\in{\cal L}_1
\end{array} \right. $$

$${\cal B}=\frac{1}{2\pi^2}\sum_{n'\in {\cal L}_2}\sqrt{\xi_{i'}\xi_{j'}}z_{n'} z_{m'}$$
$$\bar{\cal B}=\frac{1}{2\pi^2}\sum_{n'\in {\cal L}_2} \sqrt{\xi_{i'}\xi_{j'}}\bar z_{n'} \bar z_{m'}$$

For the notational simplicity, $I,\theta,H$ refer to $I_+,\theta_+,H\circ\Psi\circ\Phi$. Where  $P$ is just $G$ with the $(q_{i_1}, \cdots, q_{i_b}, \bar
q_{i_1}, \cdots, \bar q_{i_b}, q_n, \bar q_n )$-variables
expressed in terms of the $(\theta,I, z_n,\bar z_n)$ variables.
\begin{eqnarray}
H=&&\la \omega,I\ra+\sum_{[n]} \la A_{[n]} z_{[n]}, \bar z_{[n]}\ra\label{Am1}\\
&+&\sum_{{n'\in {\cal L}_2}}[ (\Omega_{n'}-\omega_{i'})z_{n'}\bar
z_{n'}+(\Omega_{m'}-\omega_{j'})z_{m'}\bar
z_{m'}]\nonumber\\
&+&\frac{1}{2\pi^2}\sum_{n'\in {\cal L}_2}\sqrt{\xi_{i'}\xi_{j'}}z_{n'} z_{m'}+\frac{1}{2\pi^2}\sum_{n'\in {\cal L}_2}\sqrt{\xi_{i'}\xi_{j'}}\bar z_{n'} \bar z_{m'}\nonumber\\
 &+&P(\theta, I, z, \bar z, \xi)\nonumber\\
 =&&N+{\cal B}+\bar{\cal B}+P\nonumber
 \end{eqnarray}

\noindent
where $A_{[n]}$ is $\sharp [n]\times\sharp [n]$ matrice in (\ref{Am1})

$$A_{[n]}=\Omega_{[n]}+(P_{i j }^{011})_{i\in [n],j\in [n]}=(\Omega_{i j }+P_{i j }^{011})_{i\in [n],j\in [n]}$$
where if $i\neq j$,$\Omega_{i j }=0$; if $i = j$,$\Omega_{i j }=\Omega_i$. When $|i-j|> K$,$P_{i j }^{011}=0$.

Next let us verify that $H=N+{\cal B}+\bar{\cal B}+P$ satisfies the assumptions
$(A1)$--$(A5)$.

\noindent
{\it  Verification of $(A1)$}:$$ \frac{\partial \omega}{\partial \xi}=\frac{1}{4\pi^2}{\left(
  \begin{array}{cccc}
1 & 2 & \cdots & 2\\
    2& 1 & \cdots & 2\\
    \vdots & \vdots & \ddots &\vdots\\
    2 & 2 & \cdots & 1
  \end{array}
\right)}_{b\times b} \\=A,
$$
It is easy to check that $\det A\neq 0$,
Thus $(A1)$ is verified.

\noindent
{\it  Verification of $(A2)$}: Take $a=3$, the proof is obvious.\\
\noindent
 {\it  Verification of $(A3)$}: This part is the same as \cite{GXY}, for the sake of completeness, we rewrite it as follows:
In the following, we only give the proof for the most complicated case.
 Let
$${\cal A}_n=A_{[n]},n\in\Z_1^2\backslash{{\cal L}_2}$$

$$ {\cal A}_n=\left(
  \begin{array}{cccc}
\Omega_{n}-\omega_{i} & -\frac{1}{2\pi^2}\sqrt{\xi_{i}\xi_{j}}\\
    \frac{1}{2\pi^2}\sqrt{\xi_{i}\xi_{j}}& -(\Omega_{m}-\omega_{j})
    \end{array}
\right) \\,n\in {\cal L}_2$$

where $(m,i,j)$ is uniquely determined by $n$. We only verify $(A3)$ for $det(\la k,\omega\ra I\pm{\cal A}_n\otimes I_2 \pm I_2\otimes {\cal A}_{n'})$ which is the most complicated. Let $A,B$ be $2\times2$ matrices, we know that $\lambda I+A\otimes I-I\otimes B=(\lambda I+A)\otimes I-I\otimes B$. Moreover,we have

\begin{Lemma}\label{matrices}
$$|A\otimes I\pm I\otimes B|={(|A|-|B|)}^2+|A|{(tr(B))}^2+|B|{(tr(A))}^2\pm(|A|+|B|)tr(A)tr(B)$$
where $|\cdot|$ denotes the determinant of the corresponding matrices.
\end{Lemma}

Case $1$. $n,n'\in {\cal L}_1$.
$$\la k,\omega\ra\pm\Omega_n\pm\Omega_{n'}$$

Set $\alpha=\varepsilon^{-3}({|i_1|^2},{|i_2|^2},\cdots,{|i_b|^2})$,$\xi=(\xi_{i_1},\xi_{i_2},\cdots,\xi_{i_b})$,$\beta=\frac{1}{4\pi^2}(2,2,\cdots,2)$,
and notice that $|n|^2+|i|^2=|m|^2+|j|^2,|n'|^2+|i'|^2=|m'|^2+|j'|^2$.

Its eigenvalues are
\begin{eqnarray*}
&&\la k,\alpha\ra\pm\varepsilon^{-3}(|n|^2+|i|^2)\pm\varepsilon^{-3}(|n'|^2+|i'|^2)+\la Ak\pm2\beta\pm2\beta,\xi\ra\\
&\pm&\frac{1}{8\pi^2}[(-\xi_i-\xi_j\pm\sqrt{{\xi_i}^2+14\xi_i\xi_j+{\xi_j}^2})\pm(-\xi_{i'}-\xi_{j'}\pm\sqrt{{\xi_{i'}}^2+14\xi_{i'}\xi_{j'}+{\xi_{j'}}^2})].
\end{eqnarray*}
If $i\neq i'$,all the eigenvalues are not identically zero due to the presence of the square root terms.\\
If $i=i'$,consequently $j=j'$,hence if the eigenvalue is
\begin{eqnarray*}
&&\la k,\alpha\ra+\varepsilon^{-3}(|n|^2+|i|^2)-\varepsilon^{-3}(|n'|^2+|i|^2)+\la Ak+2\beta-2\beta,\xi\ra\\
&+&\frac{1}{8\pi^2}[(-\xi_i-\xi_j+\sqrt{{\xi_i}^2+14\xi_i\xi_j+{\xi_j}^2})-(-\xi_{i}-\xi_{j}+\sqrt{{\xi_{i}}^2+14\xi_{i}\xi_{j}+{\xi_{j}}^2})]\\
&=&\la k,\alpha\ra+\varepsilon^{-3}(|n|^2-|n'|^2)+\la Ak,\xi\ra
\end{eqnarray*}
then $Ak\neq 0$ for $k\neq0$; if the eigenvalue is
\begin{eqnarray*}
&&\la k,\alpha\ra+\varepsilon^{-3}(|n|^2+|i|^2)+\varepsilon^{-4}(|n'|^2+|i|^2)+\la Ak+2\beta+2\beta,\xi\ra\\
&+&\frac{1}{8\pi^2}[(-\xi_i-\xi_j+\sqrt{{\xi_i}^2+14\xi_i\xi_j+{\xi_j}^2})+(-\xi_{i}-\xi_{j}-\sqrt{{\xi_{i}}^2+14\xi_{i}\xi_{j}+{\xi_{j}}^2})]\\
&=&\la k,\alpha\ra+\varepsilon^{-3}(|n|^2+|i|^2)+\varepsilon^{-3}(|n'|^2+|i|^2)+\la Ak+2\beta+2\beta,\xi\ra+\frac{1}{4\pi^2}(-\xi_i-\xi_j)\\
&=&\la k,\alpha\ra+\varepsilon^{-3}(|n|^2+|i|^2)+\varepsilon^{-3}(|n'|^2+|i|^2)+\la Ak+2\beta+2\beta+\frac{1}{4\pi^2}(-e_i-e_j),\xi\ra
\end{eqnarray*}
then when $Ak+2\beta+2\beta+\frac{1}{4\pi^2}(-e_i-e_j)=0$,all components of $k+e_i+e_j$ are equal and $(2b-1){(k+e_i+e_j)}_1+8=0(b\geq2)$,this equation has no integer solutions.Thus all eigenvalues are not identically zero.

Case $2$. $n\in {\cal L}_1,n'\in {\cal L}_2$. In this case, the eigenvalues of $(\la k,\omega\ra\pm\Omega_n) I\pm{\cal A}_{n'}$ are
 \begin{eqnarray*}
&&\la k,\alpha\ra\pm\varepsilon^{-3}(|n|^2+|i|^2)\pm\varepsilon^{-3}(|n'|^2-|i'|^2)+\la Ak\pm2\beta,\xi\ra\\
&\pm&\frac{1}{8\pi^2}[(-\xi_i-\xi_j\pm\sqrt{{\xi_i}^2+14\xi_i\xi_j+{\xi_j}^2})\pm(\xi_{i'}-\xi_{j'}\pm\sqrt{{\xi_{i'}}^2-14\xi_{i'}\xi_{j'}+{\xi_{j'}}^2})].
\end{eqnarray*}
Because $\sqrt{{\xi_{i'}}^2-14\xi_{i'}\xi_{j'}+{\xi_{j'}}^2}$ has non--zero imaginary part, there will be no small divisor.

Case $3$. $n,n'\in {\cal L}_2$. In this case, the eigenvalues of $\la k,\omega\ra I\pm{\cal A}_n\otimes I_2 \pm I_2\otimes {\cal A}_{n'}$ are
\begin{eqnarray*}
&&\la k,\alpha\ra\pm\varepsilon^{-3}(|n|^2-|i|^2)\pm\varepsilon^{-3}(|n'|^2-|i'|^2)+\la Ak,\xi\ra\\
&\pm&\frac{1}{8\pi^2}[(\xi_i-\xi_j\pm\sqrt{{\xi_i}^2-14\xi_i\xi_j+{\xi_j}^2})\pm(\xi_{i'}-\xi_{j'}\pm\sqrt{{\xi_{i'}}^2-14\xi_{i'}\xi_{j'}+{\xi_{j'}}^2})].
\end{eqnarray*}
If $i\neq i'$, all the eigenvalues are not identically zero due to the presence of the square root terms.\\
If $i=i'$, consequently $j=j'$, hence if the eigenvalue is
\begin{eqnarray*}
&&\la k,\alpha\ra+\varepsilon^{-3}(|n|^2-|i|^2)-\varepsilon^{-3}(|n'|^2-|i|^2)+\la Ak,\xi\ra\\
&+&\frac{1}{8\pi^2}[(\xi_i-\xi_j+\sqrt{{\xi_i}^2-14\xi_i\xi_j+{\xi_j}^2})-(\xi_{i}-\xi_{j}+\sqrt{{\xi_{i}}^2-14\xi_{i}\xi_{j}+{\xi_{j}}^2})]\\
&=&\la k,\alpha\ra+\varepsilon^{-3}(|n|^2-|n'|^2)+\la Ak,\xi\ra
\end{eqnarray*}
then $Ak\neq 0$ for $k\neq0$; if the eigenvalue is
\begin{eqnarray*}
&&\la k,\alpha\ra+\varepsilon^{-3}(|n|^2-|i|^2)+\varepsilon^{-3}(|n'|^2-|i|^2)+\la Ak,\xi\ra\\
&+&\frac{1}{8\pi^2}[(\xi_i-\xi_j+\sqrt{{\xi_i}^2-14\xi_i\xi_j+{\xi_j}^2})+(\xi_{i}-\xi_{j}-\sqrt{{\xi_{i}}^2-14\xi_{i}\xi_{j}+{\xi_{j}}^2})]\\
&=&\la k,\alpha\ra+\varepsilon^{-3}(|n|^2-|i|^2)+\varepsilon^{-3}(|n'|^2-|i|^2)+\la Ak,\xi\ra+\frac{1}{4\pi^2}(\xi_i-\xi_j)\\
&=&\la k,\alpha\ra+\varepsilon^{-3}(|n|^2-|i|^2)+\varepsilon^{-3}(|n'|^2-|i|^2)+\la Ak+\frac{1}{4\pi^2}(e_i-e_j),\xi\ra
\end{eqnarray*}
then when $Ak+\frac{1}{4\pi^2}(e_i-e_j)=0$, all components of $k-e_i+e_j$ are equal and $(2b-1){(k-e_i+e_j)}_1=0(b\geq2)$, this integer equation to this equation are $k=e_i-e_j$. While at this time,when $|n|\neq|m'|$,
\begin{eqnarray*}
&&\la e_i-e_j,\alpha\ra+\varepsilon^{-3}(|n|^2-|i|^2)+\varepsilon^{-3}(|n'|^2-|i|^2)\\
&=&\varepsilon^{-3}(|i|^2-|j|^2+|n|^2-|i|^2+(-|m'|^2+|j|^2))\\
&=&\varepsilon^{-3}(|n|^2-|m'|^2)\neq0
\end{eqnarray*}
Thus all the eigenvalues are not identically zero.
In other cases, the proof  is similar, so we omit it. Due to Lemma \ref{matrices}, $det(\la k,\omega\ra I\pm{\cal A}_n\otimes I_2 \pm I_2\otimes {\cal A}_{n'})$ is polynomizl function in $\xi$ of order at most four. Thus
$$|\partial^4_\xi(det(\la k,\omega\ra I\pm{\cal A}_n\otimes I_2 \pm I_2\otimes {\cal A}_{n'}))|\geq \frac12|k|\neq0.$$
By excluding some parameter set with measure $\Cal O(\gamma^{\frac14})$, we have

$$|det(\la k,\omega\ra I\pm{\cal A}_n\otimes I_2 \pm I_2\otimes {\cal A}_{n'})|\geq \frac{\gamma}{K^\tau},k\neq 0,n,n'\in{\cal L}_2,$$

Thus $(A3)$ is verified.\\

\noindent
  {\it Verification of $(A4)$}: For a given $0<r<1$ and $s=\varepsilon^{\frac 12}$, according to Lemma \ref{regularityGG}, $\|G_{\bar q}\|_\rho\leq c\|q\|_\rho^3$,
  then
  $$\sum_{n\in \Z_1^2}\|P_{w_n}\|_{\Cal O}e^{|n|\rho}+\sum_{n\in \Z_1^2}\|P_{\bar w_n}\|_{\Cal O}e^{|n|\rho}=\|P_w\|_\rho+\|P_{\bar w}\|_\rho\leq c\|q\|_\rho^3\leq
  c(|I|^{\frac 32}+\|w\|_\rho^3).$$
In addition,
$$\sup_{\|q\|_\rho<2s}\|G\|_{\Cal
O}\leq c\sup_{\|q\|_\rho<2s}\|q\|_\rho^4\le cs^4,$$ thus
$$\|P\|_{D_\rho(2r,2s),\Cal O}=\sup_{D_\rho(2r,2s)}\|P\|_{\Cal
O}\leq cs^4.$$ According to Cauchy estimates,
$$\|P_I\|_{D_\rho(r,s),\Cal O}\leq cs^2,  \|P_\theta\|_{D_\rho(r,s),\Cal O}\leq
cs^4.$$ Hence
\begin{eqnarray*} \|X_P\|_{D_\rho(r,s),\Cal O}&=&
\|P_I\|_{ D_\rho(r,s)  , \Cal O}+ \frac
1{s^2}\|P_\theta\|_{ D_\rho(r,s)  , \Cal O}\\
&+&\sup_{D_\rho(r,s)}[ \frac 1s\sum_{n\in\Z_1^2} \|P_{w_n}\|_{\cal
O}e^{|n|\rho}+ \frac 1s\sum_{n\in\Z_1^2} \|P_{\bar w_n}\|_{\cal O}
e^{|n|\rho} ]\\
&\leq&cs^2+\frac{cs^4}{s^2}+c\sup_{D_\rho(r,s)}\frac 1s(|I|^{\frac
32}+\|z\|_\rho^3)\\
&\leq&cs^2\leq c\varepsilon.
\end{eqnarray*}
   Thus $(A4)$ is verified.\\

\noindent
  {\it Verification of $(A5)$}: We only need to check $P$ satisfies $(A5)$. Recall $(\ref{P12})$. $F$ is given as
  $$F=\sum_{{{i-j+n-m=0}\atop{|i|^2-|j|^2+|n|^2-|m|^2\neq0}}\atop\sharp S\cap\{i,j,n,m\}\geq2}\frac{i}{8\pi^2(\lambda_i-\lambda_j+\lambda_n-\lambda_m)}q_i\bar q_jw_n\bar w_m.$$
  Then for $t$  large enough and $\forall c\in \Z^2\setminus\{0\}$, we have
 \begin{eqnarray*}
 &&\sum_{i,j,n,m,t}\frac{i}{8\pi^2(\lambda_i-\lambda_j+\lambda_{n+tc}-\lambda_{m+tc})}q_i\bar q_jw_{n+tc}\bar w_{m+tc}\\
 &=&\sum_{i,j,n,m,t}\frac{i}{8\pi^2(|i|^2-|j|^2+|n|^2-|m|^2+2t\la n-m,c\ra)}q_i\bar q_jw_{n+tc}\bar w_{m+tc}.
\end{eqnarray*}
 Hence,when $\la n-m,c\ra=0$,
 $$\frac{\partial^2 F}{\partial
w_{n+tc}\partial \bar w_{m+tc}}=\frac{\partial^2 F}{\partial
w_{n}\partial \bar w_{m}};$$
when $\la n-m,c\ra\neq 0$,
 $$\|\frac{\partial^2 F}{\partial
w_{n+tc}\partial \bar w_{m+tc}}-0\|\leq\frac{\varepsilon}{|t|}e^{-|n-m|\rho}.$$
 Similarly,
 $$\|\frac{\partial^2 F}{\partial
w_{n+tc}\partial  w_{m-tc}}-\lim_{t\to \infty}\frac{\partial^2 F}{\partial
w_{n+tc}\partial  w_{m-tc}}\|,\|\frac{\partial^2 F}{\partial
\bar w_{n+tc}\partial \bar w_{m-tc}}-\lim_{t\to \infty}\frac{\partial^2 F}{\partial
\bar w_{n+tc}\partial \bar w_{m-tc}}\|\leq\frac{\varepsilon}{|t|}e^{-|n+m|\rho}.$$
That is to say, $F$ satisfies T\"{o}plitz-Lipschitz property. Recalling the construction of Hamiltonian (\ref{PH}), we only need to check that
$\{G,F\}$ also satisfies the T\"{o}plitz-Lipschitz property. Lemma \ref{toplitz} in the next section shows that Poisson bracket preserves T\"{o}plitz-Lipschitz property. Thus $N+{\cal B}+\bar{\cal B}+P$ satisies $(A5)$.
Thus $(A5)$ is verified.\\
 So we have verified  all the assumptions of Theorem \ref{KAM}
  for (\ref{Hamiltoniann}).   By applying Theorem 2, we get Theorem 1.

\section{ KAM Step}

 Theorem \ref{KAM} will be proved by a KAM iteration which
involves an infinite sequence of change of variables. Each step of
KAM iteration makes the perturbation smaller than that of the previous
step at the cost of excluding a small set of parameters and
contraction of weight. We have
to prove  the convergence of the iteration  and estimate the
measure of the excluded set after infinite KAM steps.

At the $\nu$--step of the KAM iteration,  we consider Hamiltonian function
$$ H_\nu=N_\nu+{\cal B}_\nu+\bar{\cal B}_\nu+ P_\nu,$$
where $N_\nu$ is an ``integrable normal form", ${\cal B}_\nu+\bar{\cal B}_\nu+ P_\nu$ defined in $D_{\rho_\nu}(r_\nu, s_\nu)\times \Cal O_{\nu}$
with satisfying $(A1)$--$(A5)$.\\

Our goal is to construct a  map
$$\Phi_v: D_{\rho_\nu}(r_{\nu+1}, s_{\nu+1})\times\Cal O_{\nu}
 \to D_{\rho_\nu}(r_{\nu}, s_{\nu})\times\Cal O_{\nu}$$
and
 \begin{equation}\label{4.P2}
H_{\nu+1}=H_\nu\circ\Phi_\nu=N_{\nu+1}+{\cal B}_{\nu+1}+\bar{\cal B}_{\nu+1}+ P_{\nu+1}
 \end{equation}
satisfies all the above iterative assumptions $(A1)-(A5)$ on $D_{\rho_{\nu+1}}(r_{\nu+1}, s_{\nu+1})\times\Cal O_{\nu}$. Moreover,
$$
\|X_{P_{\nu+1}}\|_{D_{\rho_{\nu+1}}(r_{\nu+1}, s_{\nu+1}),
\Cal
O_{\nu}}=\|X_{H_\nu\circ\Phi_\nu}-X_{N_{\nu+1}+{\cal B}_{\nu+1}+\bar{\cal B}_{\nu+1}}\|_{D_{\rho_{\nu+1}}(r_{\nu+1},
s_{\nu+1}), \Cal O_{\nu}}\leq\varepsilon_{\nu+1}
$$

\sss To simplify notations, in what follows, the quantities
without subscripts and superscripts
refer to quantities at the $\nu^{\rm th}$ step, while the
quantities with subscript $+$ or superscript $+$ denote the corresponding quantities
at the $(\nu+1)^{\rm th}$  step. Let us
then consider the Hamiltonian \begin{eqnarray}
H=&&N+{\cal B}+\bar{\cal B}+P\nonumber\\
 =&&\la \omega,I\ra+\sum_{[n]} \la A_{[n]} z_{[n]}, \bar z_{[n]}\ra\label{Am}\\
&+&\sum_{{n'\in {\cal L}_2}}[ (\Omega_{n'}-\omega_{i'})z_{n'}\bar
z_{n'}+(\Omega_{m'}-\omega_{j'})z_{m'}\bar
z_{m'}]\nonumber\\
&+&\frac{1}{2\pi^2}\sum_{n'\in {\cal L}_2}\sqrt{\xi_{i'}\xi_{j'}}z_{n'} z_{m'}+\frac{1}{2\pi^2}\sum_{n'\in {\cal L}_2}\sqrt{\xi_{i'}\xi_{j'}}\bar z_{n'} \bar z_{m'}\nonumber\\
 &+&P(\theta, I, z, \bar z, \xi)\nonumber
 \end{eqnarray}\noindent defined in $D_\rho(r, s)\times\Cal O$.
We assume that  $|k|\le K$,

\[ |\langle k,\omega\rangle|\ge \frac{\gamma}{K^\tau}, k\neq 0\]
 \[|\langle k,\omega\rangle +\widetilde{\lambda}_j|\ge \frac{\gamma}{K^\tau},j\in{[n]} \]
 \[|\langle k,\omega\rangle \pm \widetilde{\lambda}_i\pm \widetilde{\lambda}_j|\ge \frac{\gamma}{K^\tau},i\in{[m]},j\in{[n]}\]
 where $\widetilde{\lambda}_i,\widetilde{\lambda}_j$ are eigenvalues.
  \noindent

$$|det(\la k,\omega\ra I\pm{\cal A}_n\otimes I_2 \pm I_2\otimes {\cal A}_{n'})|\geq \frac{\gamma}{K^\tau},k\neq 0,n,n'\in{\cal L}_2,$$

where

$$ {\cal A}_n=\left(
  \begin{array}{cccc}
\Omega_{n}-\omega_{i} & -\frac{1}{2\pi^2}\sqrt{\xi_{i}\xi_{j}}\\
    \frac{1}{2\pi^2}\sqrt{\xi_{i}\xi_{j}}& -(\Omega_{m}-\omega_{j})
    \end{array}
\right) \\,n\in {\cal L}_2$$
where $(n,m)$ are resonant pairs, $(i,j)$
are uniquely determined by $(n,m)$ in ${\cal L}_2$.\\
   Moreover, $N+{\cal B}+\bar{\cal B}+P$ satisfies $(A4),(A5)$.

 \sss \noindent {\bf Remark } The assumption $(A5)$ makes the  measure estimate available at each
KAM step.

  We now let $0<r_+<r$ and define
 \begin{equation}\label{4.51}
 s_+=\frac 14s\varepsilon^{\frac 13}, \quad
 \varepsilon_+=c\gamma^{-5}K^{5(\tau+1)}(r-r_+)^{-c}
\varepsilon^{\frac {4}{3}}.
\end{equation}
 Here and later, the letter $c$  denotes suitable (possibly
different) constants that do not depend on the iteration steps.

We now describe how  to construct a set $\Cal O_+\subset \Cal O$
and a change of variables
 $\Phi: D_+\times\Cal O_+=D_\rho(r_+, s_+)\times \Cal O_+\to  D_\rho(r,s)  \times \Cal O$
 such that the transformed
Hamiltonian $H_+=N_++{\cal B}_++\bar{\cal B}_++P_+\equiv H\circ \Phi$ satisfies all the
above iterative assumptions with new parameters $s_+,
\varepsilon_+, r_+$ and with $\xi\in \Cal O_+$.

\subsection{Solving the linearized equations}\label{4.1}

  Expand $P$ into the  Fourier-Taylor series
$$P=\sum_{k,l,\alpha,\beta} P_{kl\alpha\beta}\kth I^lw^\alpha\bar w^\beta$$
where $k\in \Z^b, l\in \N^b$ and the multi--indices $\alpha$ and
$\beta $ run over the set of all infinite dimensional vectors
$\alpha\equiv (\cdots,\alpha_n,\cdots)_{n\in\Z_1^2}$, $\beta
\equiv(\cdots, \beta _n, \cdots)_{n\in\Z_1^2}$ with finitely many
nonzero components of positive integers.

Let $R$ be the truncation of $P$ given by
$$R(\theta,I,z,\bar z)=R_0+R_1+R_2$$
where
$$R_0=\sum_{|k|\le K,|l|\le 1} P_{kl00}\kth I^l$$
\begin{eqnarray}
R_1&=&\sum_{{|k|\le K,n'\in {\cal L}_2}} ({P}_{n'}^{k10}z_{n'}+{P}_{m'}^{k10} z_{m'}+{P}^{k01}_{n'}\bar z_{n'}+{P}^{k01}_{m'}\bar z_{m'}) \kth
\nonumber\\
&+&\sum_{{|k|\le K,[n]}} (\la R_{[n]}^{k10}, z_{[n]}\ra+\la R^{k01}_{[n]},\bar z_{[n]} \ra) \kth
\nonumber
\end{eqnarray}

\begin{eqnarray}
R_2&=&\sum_{|k|\le K,n\in {\cal L}_2,n'\in {\cal L}_2}({P}^{k11}_{nn'}z_n \bar{z}_{n'}+
{P}^{k11}_{mn'}z_m \bar{z}_{n'}+{P}^{k11}_{nm'}z_n \bar{z}_{m'}+{P}^{k11}_{mm'}z_m\bar{z}_{m'}\nonumber\\
&&
+{P}^{k11}_{n' n}z_{n'} \bar{z}_n
+{P}^{k11}_{m' n}z_{m'} \bar{z}_n+{P}^{k11}_{n' m}z_{n'} \bar{z}_m+{P}^{k11}_{m' m}z_{m'} \bar{z}_m)\kth\nonumber\\
&+&\sum_{|k|\le K,n\in {\cal L}_2,n'\in {\cal L}_2} ({P}^{k20}_{nn'}z_{n}
z_{n'}+{P}^{k20}_{mn'}z_{m}
z_{n'}+{P}^{k20}_{nm'}z_{n}
z_{m'}+{P}^{k20}_{mm'}z_{m}
z_{m'})\kth\nonumber\\
&+&\sum_{|k|\le K,n\in {\cal L}_2,n'\in {\cal L}_2} ({P}^{k02}_{nn'}\bar{z}_{n}
\bar{z}_{n'}+{P}^{k02}_{mn'}\bar{z}_{m}
\bar{z}_{n'}+{P}^{k02}_{nm'}\bar{z}_{n}
\bar{z}_{m'}+{P}^{k02}_{mm'}\bar{z}_{m}
\bar{z}_{m'})\kth\nonumber\\
&+&\sum_{|k|\le K,[n],[m]}(\la R^{k20}_{[m][n]}z_{[n]}, z_{[m]}\ra+\la R^{k02}_{[m][n]}\bar z_{[n]},\bar z_{[m]}\ra)\kth\nonumber\\
&+&\sum_{|k|\le K,[n],[m]}\la R^{k11}_{[m][n]}z_{[n]},\bar
z_{[m]}\ra\kth\nonumber\\
&+&\sum_{|k|\le K,n\in \Z^2_1\setminus{\cal L}_2,n'\in {\cal L}_2}({P}^{k20}_{nn'}z_nz_{n'}+{P}^{k20}_{nm'}z_nz_{m'})\kth\nonumber\\
&+&\sum_{|k|\le K,n\in \Z^2_1\setminus{\cal L}_2,n'\in {\cal L}_2}({P}^{k02}_{nn'}\bar{z}_n\bar{z}_{n'}+{P}^{k0 2}_{nm'}\bar{z}_n\bar{z}_{m'})\kth\nonumber\\
&+&\sum_{|k|\le K,n\in \Z^2_1\setminus{\cal L}_2,n'\in {\cal L}_2}({P}^{k11}_{nn'}z_n\bar{z}_{n'}+{P}^{k11}_{nm'}z_n\bar{z}_{m'}+{P}^{k11}_{n'n}z_{n'}\bar{z}_n+{P}^{k11}_{m'n}z_{m'}\bar{z}_n)\kth\nonumber
\end{eqnarray}
where $P_{n}^{k10}=P_{kl\alpha\beta}$ with $\alpha=e_n, \beta=0$, here $e_n$ denotes the vector with the $n^{\rm th}$ component
 being $1$ and the other components being zero;
 $P_{n}^{k01}=P_{kl\alpha\beta}$ with $\alpha=0, \beta=e_n$;
 $P^{k20}_{nm}=P_{kl\alpha\beta}$ with $\alpha=e_n+e_m, \beta=0$;
$P^{k11}_{nm}=P_{kl\alpha\beta}$ with $\alpha=e_n, \beta=e_m$;
$P^{k02}_{nm}=P_{kl\alpha\beta}$ with $\alpha=0, \beta=e_n+e_m.$

Where,$R_{[n]}^{k10}$,$R_{[n]}^{k01}$,$R_{[m][n]}^{k20}$,$R_{[m][n]}^{k02}$ and $R_{[m][n]}^{k11}$ are, respectively, $\sharp[n]\times 1,\sharp[n]\times 1,\sharp[m]\times\sharp[n],\sharp[m]\times\sharp[n],\sharp[m]\times\sharp[n]$ matrices

$$R_{[n]}^{k10}=(P_{i}^{k10})_{i\in[n]},R_{[n]}^{k01}=(P_{i}^{k01})_{i\in[n]},|[n]|\leq K,$$
$$R_{[m][n]}^{k20}=(R_{ij}^{k20})_{i\in[m],j\in[n]}$$where if $|i+j|\leq K$, $R_{ij}^{k20}=P_{ij}^{k20}$; if $|i+j|> K$, $R_{ij}^{k20}=0$,
$$R_{[m][n]}^{k02}=(R_{ij}^{k02})_{i\in[m],j\in[n]}$$where if $|i+j|\leq K$ ,$R_{ij}^{k02}=P_{ij}^{k02}$; if $|i+j|> K$ ,$R_{ij}^{k02}=0$,
$$R_{[m][n]}^{k11}=(R_{ij}^{k11})_{i\in[m],j\in[n]}$$where if $|i-j|\leq K$ ,$R_{ij}^{k11}=P_{ij}^{k11}$; if $|i-j|> K$, $R_{ij}^{k11}=0$.

 Rewrite $H$ as $ H=N+{\cal B}+\bar{\cal B}+R+(P-R)$. By the choice of $s_+$ in
(\ref{4.51}) and the definition of the norms, it follows
immediately that
\begin{equation}
\label{4.9} \|X_R\|_{ D_\rho(r,s)  ,\Cal O}\le \| X_P\|_{ D_\rho(r,s) ,\Cal
O}\le \varepsilon.
\end{equation}

for any $\frac{r_0}{2}<\rho\leq r$. In the next, we prove that for $\frac{r_0}{2}<\rho\leq r_+$
$$\|H_{(P-R)}\|_{ D_\rho(r_+,s)  ,\Cal O}<c\varepsilon_+$$
In fact, $P-R=P^*+h.o.t.$, where
\begin{eqnarray*}P^*&=&\sum_{|n|>K}[P^{k10}_n(\theta)w_n+P^{k01}_n(\theta)\bar w_n]\\
&+&\sum_{|n+m|>K}[P^{k20}_{nm}(\theta)w_nw_m+P^{k02}_{nm}(\theta)\bar w_n\bar w_m]+\sum_{|n-m|>K}P^{k11}_{nm}(\theta)w_n\bar w_m
\end{eqnarray*} be the linear and quadratic terms in the perturbation. By virtue of (\ref{4.51}), the decay property of $P$, $\|X_P\|_{ D_\rho(r,s)  ,\Cal O}\leq\varepsilon$, and Cauchy estimates, one has that for $\rho\leq r_+$

\begin{eqnarray*}
&&\|X_{P_*}\|_{ D_\rho(r_+,s)  ,\Cal O}\\
&\leq &(r-r_+)^{-1}(\sum_{|n|>K}\varepsilon e^{-|n|r}e^{|n|\rho}+\sum_{|n+m|>K}\varepsilon e^{-|n+m|r}|\bar w_m|e^{|n+m|\rho}+\sum_{|n-m|>K}\varepsilon e^{-|n-m|r}|w_m|e^{|n-m|\rho})\\
&\leq &(r-r_+)^{-1}(\sum_{|n|>K}\varepsilon e^{-|n|r}e^{|n|\rho}+\sum_{|n|>K,m}\varepsilon e^{-|n|r}|w_m|e^{|n|\rho}e^{m\rho})\\
&\leq &(r-r_+)^{-1}\sum_{|n|>K}\varepsilon e^{-|n|(r-\rho)}\\
&\leq &(r-r_+)^{-1}\varepsilon e^{-K(r-\rho)}\\
&\leq &\varepsilon_+\\
\end{eqnarray*}

Moreover, we take $s_+\ll s$ such that in a domain $D_\rho(r, s_+)$,
\begin{equation} \label{4.10} \| X_{(P-R)}\|_ {D_\rho(r, s_+)}
\lep \varepsilon_+.
\end{equation}

\sss In the following, we will look for  an $F$, defined in a domain $D_+=D_\rho(r_+, s_+)$, such that the time one
map $\phi^1_F$ of the Hamiltonian vector field  $X_F$ defines a
map from $D_+\to D$ and transforms $H$ into $H_+$. More precisely,
by second order Taylor formula, we have
\begin{eqnarray}
H\circ \phi^1_F &=&(N+{\cal B}+\bar{\cal B}+ R)\circ \phi_F^1+(P-R)\circ \phi^1_F\nonumber\\
&=& N+{\cal B}+\bar{\cal B}+ \{N+{\cal B}+\bar{\cal B},F\}+R\nonumber\\
&+&\int_0^1 (1-t)\{\{N+{\cal B}+\bar{\cal B},F\},F\}\circ \phi_F^{t}dt\nonumber\\
&+&\int_0^1
\{R,F\}\circ \phi_F^{t}dt +(P-R)\circ \phi^1_F
\label{4.11}\\
&=& N_++{\cal B}_++\bar{\cal B}_++P_+ +\{N+{\cal B}+\bar{\cal B},F\}+R\nonumber\\
& -&P_{0000}-\la\hat{\omega},
I\ra-\sum_{[n]}\la P_{[n][n]}^{011}z_{[n]},\bar z_{[n]}\ra-\sum_{n'\in {\cal L}_2}({P}_{n'n'}^{011}z_{n'}\bar z_{n'}+{P}_{m'm'}^{011}z_{m'}\bar z_{m'})-\hat{{\cal B}}-\hat{\bar{\cal B}},\nonumber
\end{eqnarray} where
\[
\hat{\omega}= \int\frac{\partial P}{\partial I}d\theta|_{ z=\bar z= 0,
I=0},
\]

$$\hat{{\cal B}}=\sum_{n'\in {\cal L}_2}{P}^{020}_{n'm'}z_{n'} z_{m'}$$
$$\hat{\bar{\cal B}}=\sum_{n'\in {\cal L}_2} {P}^{002}_{n'm'}\bar z_{n'} \bar z_{m'}$$

 \beq\label{N_+}
 N_+= N+P_{0000}+\la\hat{\omega}, I\ra+ \sum_{[n]}\la P_{[n][n]}^{011}z_{[n]},\bar z_{[n]}\ra+\sum_{n'\in {\cal L}_2}({P}_{n'n'}^{011}z_{n'}\bar z_{n'}+{P}_{m'm'}^{011}z_{m'}\bar z_{m'}),
\eeq

\beq\label{B+}
{\cal B}_+= {\cal B}+\hat{{\cal B}},
\eeq  \beq\label{B1+}
\bar{{\cal B}}_+= \bar{{\cal B}}+\hat{\bar{{\cal B}}}= \bar{{\cal B}}+\bar{\hat{{\cal B}}},
\eeq
\beq \label{P_+} P_+=\int_0^1 (1-t)\{\{N+{\cal B}+\bar{\cal B},F\},F\}\circ
\phi_F^{t}dt+\int_0^1 \{R,F\}\circ \phi_F^{t}dt +(P-R)\circ
\phi^1_F. \eeq

 We shall
find a function $F$:

$$F(\theta, I, z,\bar z)=F_0+F_1+F_2$$

where
$$F_0=\sum_{0<|k|\le K,|l|\le 1} F_{kl00}\kth I^l$$

\begin{eqnarray}
F_1&=&\sum_{{|k|\le K,n'\in {\cal L}_2}} ({F}_{n'}^{k10}z_{n'}+{F}_{m'}^{k10} z_{m'}+{F}^{k01}_{n'}\bar z_{n'}+{F}^{k01}_{m'}\bar z_{m'}) \kth
\nonumber\\
&+&\sum_{{|k|\le K,[n]}} (\la F_{[n]}^{k10}, z_{[n]}\ra+\la F^{k01}_{[n]},\bar z_{[n]} \ra) \kth
\nonumber
\end{eqnarray}

\begin{eqnarray}
F_2&=&\sum_{|k|\le K,n\in {\cal L}_2,n'\in {\cal L}_2,|k|+|n-n'|\neq0} ({F}^{k11}_{n'n}z_{n'}
\bar{z}_{n}
+{F}^{k11}_{nn'}z_{n}\bar{z}_{n'})\kth\nonumber\\
&+&\sum_{|k|\le K,n\in {\cal L}_2,n'\in {\cal L}_2,|k|+|n-m'|\neq0} ({F}^{k11}_{m'n}z_{m'}\bar{z}_{n}+{F}^{k11}_{nm'}z_{n}\bar{z}_{m'}
)\kth\nonumber\\
&+&\sum_{|k|\le K,n\in {\cal L}_2,n'\in {\cal L}_2,|k|+|m-n'|\neq0} ({F}^{k11}_{n'm}z_{n'}\bar{z}_{m}+{F}^{k11}_{mn'}z_{m}\bar{z}_{n'}
)\kth\nonumber\\
&+&\sum_{|k|\le K,n\in {\cal L}_2,n'\in {\cal L}_2,|k|+|m-m'|\neq0} ({F}^{k11}_{m'm}z_{m'}\bar{z}_{m}
+{F}^{k11}_{mm'}z_{m}\bar{z}_{m'}
)\kth\nonumber\\
&+&\sum_{|k|\le K,n\in {\cal L}_2,n'\in {\cal L}_2,|k|+|n-m'|\neq0(or)|k|+|n'-m|\neq0} ({F}^{k20}_{n'n}z_{n'}
z_{n}+{F}^{k02}_{n'n}\bar{z}_{n'}\bar{z}_{n})\kth\nonumber\\
&+&\sum_{|k|\le K,n\in {\cal L}_2,n'\in {\cal L}_2,|k|+|m-m'|\neq0(or)|k|+|n'-n|\neq0} ({F}^{k20}_{m'n}z_{m'}z_{n}+
{F}^{k02}_{m'n}\bar{z}_{m'}\bar{z}_{n})\kth\nonumber\\
&+&\sum_{|k|\le K,n\in {\cal L}_2,n'\in {\cal L}_2,|k|+|m-m'|\neq0(or)|k|+|n'-n|\neq0} ({F}^{k20}_{n'm}z_{n'}z_{m}+{F}^{k02}_{n'm}\bar{z}_{n'}\bar{z}_{m})\kth\nonumber\\
&+&\sum_{|k|\le K,n\in {\cal L}_2,n'\in {\cal L}_2,|k|+|n-m'|\neq0(or)|k|+|n'-m|\neq0} ({F}^{k20}_{m'm}z_{m'}z_{m}+{F}^{k02}_{m'm}\bar{z}_{m'}\bar{z}_{m})\kth\nonumber\\
&+&\sum_{|k|\le K,[n],[m]}(\la F^{k20}_{[m][n]}z_{[n]}, z_{[m]}\ra+\la F^{k02}_{[m][n]}\bar z_{[n]},\bar z_{[m]}\ra)\kth\nonumber\\
&+&\sum_{|k|\le K,[n],[m],|k|+ ||n|-|m|| \neq 0} \la F^{k11}_{[m][n]}z_{[n]},\bar
z_{[m]}\ra\kth\nonumber\\
&+&\sum_{|k|\le K,n\in \Z^2_1\setminus{\cal L}_2,n'\in {\cal L}_2}({F}^{k20}_{nn'}z_nz_{n'}+{F}^{k20}_{nm'}z_nz_{m'})\kth\nonumber\\
&+&\sum_{|k|\le K,n\in \Z^2_1\setminus{\cal L}_2,n'\in {\cal L}_2}({F}^{k02}_{nn'}\bar{z}_n\bar{z}_{n'}+{F}^{k0 2}_{nm'}\bar{z}_n\bar{z}_{m'})\kth\nonumber\\
&+&\sum_{|k|\le K,n\in \Z^2_1\setminus{\cal L}_2,n'\in {\cal L}_2}({F}^{k11}_{nn'}z_n\bar{z}_{n'}+{F}^{k11}_{nm'}z_n\bar{z}_{m'}+{F}^{k11}_{n'n}z_{n'}\bar{z}_n+{F}^{k11}_{m'n}z_{m'}\bar{z}_n)\kth\nonumber
\end{eqnarray}

\noindent
where $F_{[n]}^{k10}$, $F_{[n]}^{k01}$, $F_{[m][n]}^{k20}$, $F_{[m][n]}^{k02}$ and $F_{[m][n]}^{k11}$ are, respectively, $\sharp[n]\times 1, \sharp[n]\times 1, \sharp[m]\times\sharp[n], \sharp[m]\times\sharp[n], \sharp[m]\times\sharp[n]$ matrices

$$F_{[n]}^{k10}=(F_{i}^{k10})_{i\in[n]},F_{[n]}^{k01}=(F_{i}^{k01})_{i\in[n]},|[n]|\leq K,$$
$$F_{[m][n]}^{k20}=(f_{ij}^{k20})_{i\in[m],j\in[n]}$$ where if $|i+j|\leq K$, $f_{ij}^{k20}=F_{ij}^{k20}$; if $|i+j|> K$,$f_{ij}^{k20}=0$,
$$F_{[m][n]}^{k02}=(f_{ij}^{k02})_{i\in[m],j\in[n]}$$ where if $|i+j|\leq K$, $f_{ij}^{k02}=F_{ij}^{k02}$; if $|i+j|> K$,$f_{ij}^{k02}=0$,
$$F_{[m][n]}^{k11}=(f_{ij}^{k11})_{i\in[m],j\in[n]}$$ where if $|i-j|\leq K$, $f_{ij}^{k11}=F_{ij}^{k11}$; if $|i-j|> K$,$f_{ij}^{k11}=0$.

 satisfying the equation
\begin{equation}\label{4.13}\begin{array}{rlcl}
&\{N+{\cal B}+\bar{\cal B},F\}+R-P_{0000}-\la\hat{\omega},
I\ra-\displaystyle\sum_{[n]}\la P_{[n][n]}^{011}z_{[n]},\bar z_{[n]}\ra-\hat{{\cal B}}-\hat{\bar{\cal B}}\\&-\displaystyle\sum_{n'\in {\cal L}_2}({P}_{n'n'}^{011}z_{n'}\bar z_{n'}+{P}_{m'm'}^{011}z_{m'}\bar z_{m'})=0.
\end{array}
\end{equation}

\begin{Lemma}\label{Lem4.01} $F$ satisfies  (\ref{4.13})  if the Fourier coefficients of $F_0,F_1$ are defined
by the following equations \beq\label{4.14}\begin{array}{rlcl}
 &(\la k,\omega\ra )F_{kl00}&=&
{\rm i} P_{kl00}, \quad |l|\le 1,0<|k|\le K,
\\
&(\la k,\omega\ra I -  A_{[n]})F^{k10}_{[n]}&=&{\rm i} P^{k10}_{[n]},\quad |k|\le K,n\in \Z^2_1\setminus{\cal L}_2,\\
&(\la k,\omega\ra I +  A_{[n]})F^{k01}_{[n]}&=&{\rm i} R^{k01}_{[n]},\quad |k|\le K,n\in \Z^2_1\setminus{\cal L}_2,\\
&(\la k,\omega\ra I-{\cal A}_{n'})({F}_{n'}^{k10},{F}_{m'}^{k01})^T&=&i({P}_{n'}^{k10},{P}_{m'}^{k01})^T,|k|\le K,n'\in {\cal L}_2,\\
&(\la k,\omega\ra I+{\cal A}_{n'})({F}_{n'}^{k01},{F}_{m'}^{k10})^T&=&i({P}_{n'}^{k01},{P}_{m'}^{k10})^T,|k|\le K,n'\in {\cal L}_2.
\end{array}\eeq
\end{Lemma}

The Fourier coefficients of $F_2$ are defined by the following Lemmas\\
{\bf Case 1: $n,m\in \Z^2_1\setminus{\cal L}_2$}
\begin{Lemma}\label{Lem4.02} $F$ satisfies  (\ref{4.13})  if the Fourier coefficients of $F_2$ are defined
by the following equations \beq\label{4.14}\begin{array}{rlcl}
&(\la k,\omega\ra I - A_{[m]})F^{k20}_{[m][n]}-F^{k20}_{[m][n]}A_{[n]}&=&{\rm i}
R^{k20}_{[m][n]},\\
&(\la k,\omega\ra I - A_{[m]})F^{k11}_{[m][n]}+F^{k11}_{[m][n]}A_{[n]}&=&{\rm i}
R^{k11}_{[m][n]},\quad |k|+ ||n|-|m|| \neq 0,\\
&(\la k,\omega\ra I + A_{[m]})F^{k02}_{[m][n]}+F^{k02}_{[m][n]}A_{[n]}&=&{\rm i}
R^{k02}_{[m][n]}.\\
\end{array}\eeq
\end{Lemma}
{\bf Case 2: $n\in \Z^2_1\setminus{\cal L}_2,n'\in {\cal L}_2$}
\begin{Lemma}\label{Lem4.03} $F$ satisfies  (\ref{4.13})  if the Fourier coefficients of $F_2$ are defined
by the following equations \beq\label{4.14}\begin{array}{rlcl}
&[(\la k,\omega\ra-\Omega_n) I-{\cal A}_{n'}]({F}_{nn'}^{k20},{F}_{nm'}^{k11})^T&=&i({P}_{nn'}^{k20},{P}_{nm'}^{k11})^T,\\
&[(\la k,\omega\ra+\Omega_n) I+{\cal A}_{n'}]({F}_{nn'}^{k02},{F}_{m'n}^{k11})^T&=&i({P}_{nn'}^{k02},{P}_{m'n}^{k11})^T,\\
&[(\la k,\omega\ra-\Omega_n) I+{\cal A}_{n'}]({F}_{nn'}^{k11},{F}_{nm'}^{k20})^T&=&i({P}_{nn'}^{k11},{P}_{nm'}^{k20})^T,\\
&[(\la k,\omega\ra+\Omega_n) I-{\cal A}_{n'}]({F}_{n'n}^{k11},{F}_{m'n}^{k02})^T&=&i({P}_{n'n}^{k11},{P}_{m'n}^{k02})^T.
\end{array}\eeq
\end{Lemma}
{\bf Case 3: $n,n'\in {\cal L}_2$}

\begin{Lemma}\label{Lem4.04} $F$ satisfies  (\ref{4.13})  if the Fourier coefficients of $F_2$ are defined
by the following equations

\[
\begin{array}{rlcl}
&(\la k,\omega\ra  I-{\cal A}_{n}\otimes I+I\otimes{\cal A}_{n'})({F}_{n n'}^{k11},{F}_{n m'}^{k20},{F}_{m n '}^{k02},{F}_{m'm }^{k11})^T&=&i({P}_{n n'}^{k11},{P}_{n m'}^{k20},{P}_{m n '}^{k02},{P}_{m 'm}^{k11})^T,\\

&(\la k,\omega\ra  I+{\cal A}_{n}\otimes I-I\otimes{\cal A}_{n'})({F}_{n'n}^{k11},{F}_{m'n}^{k02},{F}_{n 'm}^{k20},{F}_{m m' }^{k11})^T&=&i({P}_{n'n}^{k11},{P}_{m'n}^{k02},{P}_{n 'm}^{k20},{P}_{mm'}^{k11})^T,\\

&(\la k,\omega\ra  I-{\cal A}_{n}\otimes I-I\otimes{\cal A}_{n'})({F}_{n n'}^{k20},{F}_{n m'}^{k11},{F}_{n 'm}^{k11},{F}_{m m '}^{k02})^T&=&i({P}_{n n'}^{k20},{P}_{n m'}^{k11},{P}_{n 'm}^{k11},{P}_{m m '}^{k02})^T,\\

&(\la k,\omega\ra  I+{\cal A}_{n}\otimes I+I\otimes{\cal A}_{n'})({F}_{n n'}^{k02},{F}_{m'n}^{k11},{F}_{m n '}^{k11},{F}_{m' m }^{k20})^T&=&i({P}_{n n'}^{k02},{P}_{m'n}^{k11},{P}_{m n '}^{k11},{P}_{m ' m}^{k20})^T.
\end{array}
\]
\end{Lemma}

In the following, we only give the proof for the most complicated case.
\proof
Inserting $F$ into (\ref{4.13}).By comparing the Fourier coefficients,more precisely,\\
if $({n'},{m'})$ is a resonant pair in ${\cal L}_2$, we have
$$\sum_{\quad |k|\le K,n'\in {\cal L}_2}[\la k,\omega\ra-(\Omega_{n'}-\omega_{i'})]{F}^{k10}_{n'}z_{n'}\kth-\frac{1}{2\pi^2}\sqrt{\xi_{i'}\xi_{j'}}{F}^{k01}_{m'}z_{n'}\kth=i\sum_{\quad |k|\le K,{n'}\in {\cal L}_2} {P}^{k10}_{n'}z_{n'}\kth$$
$$\sum_{\quad |k|\le K,{n'}\in {\cal L}_2}[\la k,\omega\ra+(\Omega_{m'}-\omega_{j'})]{F}^{k01}_{m'}z_{m'}\kth-\frac{1}{2\pi^2}\sqrt{\xi_{i'}\xi_{j'}}F^{k10}_{n'}z_{m'}\kth=i\sum_{\quad |k|\le K,{n'}\in {\cal L}_2} P^{k01}_{m'}z_{m'}\kth$$
we rewrite in matrix form
$$(\la k,\omega\ra I+{\cal A}_{n'})({F}_{n'}^{k01},{F}_{m'}^{k10})^T=i({P}_{n'}^{k01},{P}_{m'}^{k10})^T,|k|\le K,n'\in {\cal L}_2,$$
similarly,form
$$(\la k,\omega\ra I-{\cal A}_{n'})({F}_{n'}^{k10},{F}_{m'}^{k01})^T=i({P}_{n'}^{k10},{P}_{m'}^{k01})^T,|k|\le K,n'\in {\cal L}_2,$$
If $({n},{m})$ and $({n'},{m'})$ are resonant pairs in ${\cal L}_2$, comparing the Fourier cofficients, we have that $({F}_{n n'}^{k11},{F}_{n m'}^{k20},{F}_{m n '}^{k02},{F}_{m'm }^{k11})^T$ satisfy

\begin{eqnarray*}
&&[\la k,\omega\ra-(\Omega_{n}-\omega_{i})+(\Omega_{n'}-\omega_{i'})]{F}^{k11}_{nn'}\kth-\frac{1}{2\pi^2}\sqrt{\xi_{i'}\xi_{j'}}{F}^{k20}_{nm'}
\kth+\frac{1}{2\pi^2}\sqrt{\xi_{i}\xi_{j}}{F}^{k02}_{mn'}\kth\\&&=i {P}^{k11}_{nn'}\kth
\end{eqnarray*}

similarly,
\begin{eqnarray*}
&&[\la k,\omega\ra-(\Omega_{n}-\omega_{i})-(\Omega_{m'}-\omega_{j'})]{F}^{k20}_{nm'}\kth+\frac{1}{2\pi^2}\sqrt{\xi_{i'}\xi_{j'}}{F}^{k11}_{nn'}
\kth+\frac{1}{2\pi^2}\sqrt{\xi_{i}\xi_{j}}{F}^{k11}_{m'm}\kth\\&&=i {P}^{k20}_{nm'}\kth
\end{eqnarray*}
\begin{eqnarray*}
&&[\la k,\omega\ra+(\Omega_{m}-\omega_{j})+(\Omega_{n'}-\omega_{i'})]{F}^{k02}_{m n'}\kth-\frac{1}{2\pi^2}\sqrt{\xi_{i'}\xi_{j'}}{F}^{k11}_{m'm}
\kth-\frac{1}{2\pi^2}\sqrt{\xi_{i}\xi_{j}}{F}^{k11}_{nn'}\kth\\&&=i {P}^{k02}_{m n'}\kth
\end{eqnarray*}
\begin{eqnarray*}
&&[\la k,\omega\ra+(\Omega_{m}-\omega_{j})-(\Omega_{m'}-\omega_{j'})]{F}^{k11}_{m' m}\kth+\frac{1}{2\pi^2}\sqrt{\xi_{i'}\xi_{j'}}{F}^{k02}_{m n'}
\kth+\frac{1}{2\pi^2}\sqrt{\xi_{i}\xi_{j}}{F}^{k20}_{nm'}\kth\\&&=i {P}^{k11}_{m' m}\kth
\end{eqnarray*}
we rewrite them into matrix form
\begin{eqnarray*}
(\la k,\omega\ra  I-{\cal A}_{n}\otimes I+I\otimes{\cal A}_{n'})({F}_{n n'}^{k11},{F}_{n m'}^{k20},{F}_{m n '}^{k02},{F}_{m'm }^{k11})^T=i({P}_{n n'}^{k11},{P}_{n m'}^{k20},{P}_{m n '}^{k02},{P}_{m 'm}^{k11})^T,
|k|\le K,n ,n'\in {\cal L}_2\\
\end{eqnarray*}
similarly, from
\begin{eqnarray*}
(\la k,\omega\ra  I+{\cal A}_{n}\otimes I-I\otimes{\cal A}_{n'})({F}_{n'n}^{k11},{F}_{m'n}^{k02},{F}_{n 'm}^{k20},{F}_{m m' }^{k11})^T=i({P}_{n'n}^{k11},{P}_{m'n}^{k02},{P}_{n 'm}^{k20},{P}_{mm'}^{k11})^T,
|k|\le K,n ,n'\in {\cal L}_2\\
(\la k,\omega\ra  I-{\cal A}_{n}\otimes I-I\otimes{\cal A}_{n'})({F}_{n n'}^{k20},{F}_{n m'}^{k11},{F}_{n 'm}^{k11},{F}_{m m '}^{k02})^T=i({P}_{n n'}^{k20},{P}_{n m'}^{k11},{P}_{n 'm}^{k11},{P}_{m m '}^{k02})^T,
|k|\le K,n ,n'\in {\cal L}_2\\
(\la k,\omega\ra  I+{\cal A}_{n}\otimes I+I\otimes{\cal A}_{n'})({F}_{n n'}^{k02},{F}_{m'n}^{k11},{F}_{m n '}^{k11},{F}_{m' m }^{k20})^T=i({P}_{n n'}^{k02},{P}_{m'n}^{k11},{P}_{m n '}^{k11},{P}_{m ' m}^{k20})^T,
|k|\le K,n ,n'\in {\cal L}_2\\
\end{eqnarray*}
In other cases, the proof  is similar, so we omit it. Thus these Lemmas are obtained.
\qed

\noindent{\bf Remark.} In the case that $({n},{m})$ and $({n'},{m'})$ are resonant pairs in ${\cal L}_2$, we have that $k,({n},{m}),({n'},{m'})$ satisfy
\begin{eqnarray*}
&&\sum_{|k|\le K,n\in {\cal L}_2,n'\in {\cal L}_2,|k|+|n-n'|\neq0} ({F}^{k11}_{n'n}z_{n'}
\bar{z}_{n}+{F}^{k11}_{nn'}z_{n}\bar{z}_{n'}
)\kth\nonumber\\
&+&\sum_{|k|\le K,n\in {\cal L}_2,n'\in {\cal L}_2,|k|+|n-m'|\neq0(or)|k|+|n'-m|\neq0} ({F}^{k20}_{n'n}z_{n'}
z_{n}+{F}^{k02}_{n'm'}\bar{z}_{n'}
\bar{z}_{n})\kth\nonumber\\
&+&\cdots\cdots\nonumber
\end{eqnarray*}
\\
Consider the equations
$$Q^T_{[n]}(\la k,\omega\ra I-A_{[n]})F^{k10}_{[n]}=iQ^T_{[n]}R^{k10}_{[n]},|k|\leq K,$$
matrix $Q_{[n]}$ is the $A_{[n]}$'s orthogonal matrix
$$(\la k,\omega\ra I-Q^T_{[n]}A_{[n]}Q_{[n]})Q^T_{[n]}F^{k10}_{[n]}=iQ^T_{[n]}R^{k10}_{[n]},|k|\leq K,$$
that is
$$(\la k,\omega\ra I-\Lambda_{[n]})\hat{F}^{k10}_{[n]}=i\hat{R}^{k10}_{[n]},|k|\leq K.$$
Similarly, from
\begin{eqnarray*}
(\la k,\omega\ra I+\Lambda_{[n]})\hat{F}^{k01}_{[n]}&=&i\hat{R}^{k01}_{[n]},|k|\leq K,\\
(\la k,\omega\ra I-\Lambda_{[m]})\hat{F}^{k20}_{[m][n]}-\hat{F}^{k20}_{[m][n]}\Lambda_{[n]}&=&i\hat{R}^{k20}_{[m][n]},|k|\leq K,\\
(\la k,\omega\ra I-\Lambda_{[m]})\hat{F}^{k11}_{[m][n]}+\hat{F}^{k11}_{[m][n]}\Lambda_{[n]}&=&i\hat{R}^{k11}_{[m][n]},|k|\leq K,|k|+||n|-|m||\neq0,\\
(\la k,\omega\ra I+\Lambda_{[m]})\hat{F}^{k02}_{[m][n]}+\hat{F}^{k02}_{[m][n]}\Lambda_{[n]}&=&i\hat{R}^{k02}_{[m][n]},|k|\leq K.
\end{eqnarray*}
instead, where $A_{[n]}$ can be diagonalized by orthogonal matrix $Q_{[n]}$, that is $\Lambda_{[n]}=Q^T_{[n]}A_{[n]}Q_{[n]}$.
\begin{eqnarray*}
\hat{R}^{kx}_{[n]}&=&Q^T_{[n]}R^{kx}_{[n]},x=10,01\\
\hat{R}^{kx}_{[m][n]}&=&Q^T_{[m]}R^{kx}_{[m][n]}Q_{[n]},x=20,11,02.
\end{eqnarray*}
\begin{eqnarray*}
\hat{F}^{kx}_{[n]}&=&Q^T_{[n]}F^{kx}_{[n]},x=10,01\\
\hat{F}^{kx}_{[m][n]}&=&Q^T_{[m]}F^{kx}_{[m][n]}Q_{[n]},x=20,11,02.
\end{eqnarray*}
Now we focus on the following equations
\begin{eqnarray*}
(\la k,\omega\ra-\widetilde{\lambda}_j)\hat{F}^{k10}_{[n],j}&=&i\hat{R}^{k10}_{[n],j},|k|\leq K,j\in [n],\\
(\la k,\omega\ra+\widetilde{\lambda}_j)\hat{F}^{k01}_{[n],j}&=&i\hat{R}^{k01}_{[n],j},|k|\leq K,j\in [n],\\
(\la k,\omega\ra-\widetilde{\lambda}_i-\widetilde{\lambda}_j)\hat{F}^{k20}_{[m][n],ij}&=&i\hat{R}^{k20}_{[m][n],ij},|k|\leq K,i\in [m],j\in [n],\\
(\la k,\omega\ra-\widetilde{\lambda}_i+\widetilde{\lambda}_j)\hat{F}^{k11}_{[m][n],ij}&=&i\hat{R}^{k11}_{[m][n],ij},|k|\leq K,|k|+||n|-|m||\neq0,i\in [m],j\in [n],\\
(\la k,\omega\ra+\widetilde{\lambda}_i+\widetilde{\lambda}_j)\hat{F}^{k02}_{[m][n],ij}&=&i\hat{R}^{k02}_{[m][n],ij},|k|\leq K,i\in [m],j\in [n].
\end{eqnarray*}
In the other cases, the proof  is similar, so we omit it.
In order to solve the last three equations, we need the following elementary algebraic result from matrix theory.
\begin{Lemma}\label{Lem4.2} Let $A,B,C$ be, respectively, $n\times n,m\times m,n\times m$ matrices, and let $X$ be an $n\times m$ unknown matrix.
The matrix equation
$$AX-XB=C,$$
is solvable if and only if $I_m\otimes A-B\otimes I_n$ is nonsingular.
\end{Lemma}

For a detailed proof, we refer the reader to the Appendix in \cite{YJ}.\\
\noindent{\bf Remark.} Taking the transpose of the fourth equation in Lemma \ref{Lem4.02}, one sees that $(F^{k20}_{[m][n]})^T$ satisfies the same equation as $(F^{k20}_{[n][m]})$. Then (by the uniqueness of the solution) it follows that $(F^{k02}_{[n][m]})=(F^{k02}_{[m][n]})^T$, $(F^{-k11}_{[n][m]})=\overline{(F^{k11}_{[m][n]})^T}$

\subsection{Estimation
for coefficients of $F$}\label{4.2}
Let us consider $F^{k20}_{[m][n]}$ for instance, and the other terms can be treated in an analogous way. By the construction above, one sees that
$$F^{k20}_{[m][n],ij}=i\sum_{m_1,n_1}\frac{Q_{[m],im_1}\hat{R}^{k20}_{[m][n],m_1,n_1}Q^T_{[n],n_1j}}{\la k,\omega\ra-\widetilde{\lambda}_i-\widetilde{\lambda}_j}.$$
Then
$$|F^{k20}_{[m][n],ij}|\leq c\varepsilon\frac{K^\tau}{\gamma}e^{K^{1+3\varepsilon}\rho}e^{-\rho|i+j|}e^{-|k|r},$$
where we used the factor $e^{K^{1+3\varepsilon}\rho}$ to recover the exponential decay
under the assumption
$$K^{1+3\varepsilon}\rho=1.$$
And
$$\|F^{k20}_{[m][n]}\|\leq cK_{\nu}^{3\varepsilon}\varepsilon_{\nu+1}\frac{K_{\nu}^{5(\tau+1)}}{\gamma^{-5}}K_{\nu}^{3\varepsilon}\leq \varepsilon_{\nu+1}^{\frac13}$$
under the assumption
$$\varepsilon_{\nu+1}=c\gamma^{-5}(r_{\nu}-r_{\nu+1})^{-c}
K_{\nu}^{5(\tau+1)}\varepsilon_{\nu}^{\frac43}.$$

\subsection{Estimation on the coordinate transformation}\label{4.3}

 \sss We proceed to estimate $X_F$ and $\phi_F^1$. We
start with the following
\begin{Lemma}\label{Lem4.3}
Let $D_i=D( r_++\frac{i}4 (r-r_+), \frac i4s)$, $0 <i \le 4$. Then
\begin{equation}\label{4.20}
\|X_F\|_{D_3, \Cal O}\le c\gamma^{-5}K^{5(\tau+1)}(r-r_+)^{-c}\varepsilon.
\end{equation}
\end{Lemma}

 In the next lemma, we give some
estimates for $\phi_F^t$. The  formula (\ref{4.26}) will be used
to prove  our coordinate transformation is well defined.
Inequality
 (\ref{4.27}) will be used to check the convergence of the iteration.

\begin{Lemma}\label{Lem4.4}
Let  $\eta=\varepsilon^{\frac 13}, D_{i\eta}= D(r_++\frac
{i}4(r-r_+),\frac i4 \eta s), 0 <i \le 4$. If $\varepsilon\ll
\frac 12\gamma^{\frac{15}{2}}K^{{\frac{15}{2}}(\tau+1)}(r-r_+)^c$, we then have
\begin{equation}
\phi_F^t:  D_{2\eta} \to  D_{3\eta} ,\ \ \ -1 \le t\le 1,
\label{4.26}
\end{equation}
 Moreover,
\begin{equation}
\|D\phi_F^t-Id\|_{D_{1\eta}}< c
\gamma^{-5}K^{5(\tau+1)}(r-r_+)^{-c}\varepsilon. \label{4.27}
\end{equation}
\end{Lemma}
\proof Let
$$\|D^mF\|_{D,\Cal O}
=\max \{ \|\frac{\partial^{|i|+|l|+|\alpha|+|\beta|}}{\partial
\theta^{i}\partial I^{l}
\partial  z^\alpha\partial{\bar z^\beta}} F\|_{D, \Cal O}, |i|+|l|+|\alpha|+|\beta|=m\ge 2\}.$$

 Notice that  $F$ is a polynomial of degree 1 in $I$ and degree 2 in $z$, $\bar z$.
  From
(\ref{2.6}),
 (\ref{4.20}) and the Cauchy inequality,
 it  follows
that
 \begin{equation}\|D^mF\|_{D_2, \Cal O }<
 c \gamma^{-5}K^{5(\tau+1)}(r-r_+)^{-c}\varepsilon,\label{4.28}
\end{equation}
for any $m\ge 2$.

 To get the estimates for $\phi_F^t$, we start from the integral
equation,
$$\phi_F^t=id+\int_0^tX_F\circ \phi_F^s\,ds$$
so that $\phi_F^t: D_{2\eta} \to  D_{3\eta} ,\ \ \ -1\le t\le 1$,
which follows directly from (\ref{4.28}).
 Since
$$D\phi_F^t=Id+\int_0^t(DX_F) D\phi_F^s\,ds=
Id+\int_0^t J(D^2F) D\phi_F^s\,ds,$$ \noindent where $J$ denotes
the standard symplectic matrix $\left(\begin{array}{cc} 0&-I \\
I&0 \end{array}\right)$,  it follows that
\begin{equation}\|D\phi_F^t-Id\|\le 2\|D^2F\|<
c \gamma^{-5}K^{5(\tau+1)}(r-r_+)^{-c}\varepsilon. \label{4.29}
\end{equation}
Consequently Lemma \ref{Lem4.4} follows. \qed

\subsection{Estimation for the new normal form}\label{no4.4}

  The map $\phi_F^1$ defined above transforms $H$
into $H_+=N_++{\cal B}_++\bar{\cal B}_++P_+$(see (\ref{4.11}) and (\ref{4.13}))with
the normal form $N_+$
\begin{eqnarray*}
N_+&=& N+P_{0000}+\la\hat{\omega}, I\ra+\sum_{[n]}\la P_{[n][n]}^{011}z_{[n]},\bar z_{[n]}\ra+\sum_{n'\in {\cal L}_2}({P}_{n'n'}^{011}z_{n'}\bar z_{n'}+{P}_{m'm'}^{011}z_{m'}\bar z_{m'})
\nonumber\\
&=&\la\omega_+, I\ra
 +  \sum_{[n]}\la A_{[n]}^+z_{[n]} ,\bar z_{[n]}\ra+\sum_{n'\in {\cal L}_2}[(\Omega_{n'}^+-\omega_{i'})z_{n'}\bar z_{n'}+(\Omega_{m'}^+-\omega_{j'})z_{m'}\bar z_{m'})]
\end{eqnarray*}
 where
 \beq\label{frequenciesomega}
\omega_+=\omega+P_{0l00} (|l|=1), \eeq
$$A_{[n]}^+=A_{[n]}+R_{[n][n]}^{011}=A_{[n]}+(R_{ij}^{011})_{i\in [n],j\in [n]},|i-j|> K,R_{ij}^{011}=0;|i-j|\leq K,R_{ij}^{011}=P_{ij}^{011}$$
$$\Omega_{n'}^+=\Omega_{n'}+{P}_{n'n'}^{011},\Omega_{m'}^+=\Omega_{m'}+{P}_{m'm'}^{011},n'\in {\cal L}_2$$
Now we prove that $N_+$ shares the same properties as $N$. By
the regularity of $X_P$ and by
Cauchy estimates, then  we have
\begin{equation}\label{4.32}
|\omega_{+}-\omega|<\varepsilon, \quad
|P_{ij+}^{011}-P_{ij}^{011}|<\varepsilon e^{-|i-j|\rho}.
\end{equation}
 It follows that for $|k|\le K$,
 $$|\la k,\omega+P_{0l00}\ra|
 \ge|\la k,\omega\ra|-\varepsilon K\ge \frac{\gamma}{K^\tau}-\varepsilon K\ge \frac{\gamma}{K_+^\tau},
$$
 $$|\la k,\omega+P_{0l00}\ra  +\widetilde{\lambda}^+_j|
 \geq
|\la k,\omega\ra  +\widetilde{\lambda}_j|-\varepsilon K\geq \frac{\gamma}{K^\tau}-\varepsilon K\geq\frac{\gamma}{K_+^\tau},
$$
Similarly,we have
$$|\la k,\omega+P_{0l00}\ra  +\widetilde{\lambda}^+_i\pm\widetilde{\lambda}^+_j|
 \geq\frac{\gamma}{K_+^\tau}.$$
In other cases, the proof  is similar, so we omit it.

This means that in the
next KAM step, small denominator conditions are automatically
satisfied for $|k|\le  K$. The following bounds will be used for the measure eatimates:
$$\sup_{\xi\in \Cal O}\max_{d\leq 4}\|\partial_{\xi}^d(A_{[n]}^+-A_{[n]})\|\leq c\varepsilon$$
$$\sup_{\xi\in \Cal O}\max_{d\leq 4}|\partial_{\xi}^d(\Omega_{n'}^+-\Omega_{n'})|\leq \varepsilon$$
$$\sup_{\xi\in \Cal O}\max_{d\leq4}|\partial_\xi^d(\omega_+-\omega)|\leq \varepsilon$$
and $$|P^{011}_{ij+}-P_{ij}^{011}|_{\Cal O}\leq\varepsilon e^{-|i-j|\rho}.$$

\subsection{Estimation
for the new perturbation}\label{4.5}

Since \begin{eqnarray*} P_+&=&\int_0^1 (1-t)\{\{N+{\cal B}+\bar{\cal B},F\},F\}\circ
\phi_F^{t}dt+\int_0^1 \{R,F\}\circ \phi_F^{t}dt +(P-R)\circ
\phi^1_F\\
&=&\int_0^1 \{R(t),F\}\circ \phi_F^{t}dt +(P-R)\circ \phi^1_F,
\end{eqnarray*}
where $R(t)=(1-t)(N_++{\cal B}_++\bar{\cal B}_+-N-{\cal B}-\bar{\cal B})+tR$. Hence
$$
X_{P_+}=\int_0^1 (\phi_F^{t})^*X_{\{R(t),F\}} dt
+(\phi^1_F)^*X_{(P-R)}.
$$
 According to Lemma \ref{Lem4.4}, $$\|D\phi_F^t-Id\|_{D_{1\eta}}< c
\gamma^{-5}K^{5(\tau+1)}(r-r_+)^{-c}\varepsilon, \quad -1\le t\le 1,
$$
thus
$$\|D\phi_F^t\|_{D_{1\eta}}\le 1+\|D\phi_F^t-Id\|_{D_{1\eta}}\le 2, \quad -1\le t\le
1.$$ Due to Lemma \ref{Lem7.3},
$$ \|X_{\{R(t),F\}}\|_{D_{2\eta}}\le c
\gamma^{-5} K^{5(\tau+1)}(r-r_+)^{-c} \eta^{-2} \varepsilon^2,$$ and $$
\|X_{(P-R)}\|_{D_{2\eta}}\le c \eta \varepsilon, $$ we have  $$
\|X_{P_+}\|_{D_\rho(r_+,s_+)}\le c\eta \varepsilon + c \gamma^{-5}
K^{5(\tau+1)}(r-r_+)^{-c}\eta^{-2} \varepsilon^2\le c\varepsilon_+. $$

\subsection{Verification of $(A5)$ after one step of KAM iteration}\label{4.5}

Since
 \begin{eqnarray*}
P_+&=&P-R+\{P,F\}+\frac{1}{2!}\{\{N+{\cal B}+\bar{\cal B},F\},F\}+\frac{1}{2!}\{\{P,F\},F\}\\
&&+\cdots+ \frac{1}{n!}\{\cdots\{N+{\cal B}+\bar{\cal B},\underbrace{F\}\cdots
,F}_n\}+\frac{1}{n!}\{\cdots\{P,\underbrace{F\}\cdots
,F}_n\}+\cdots
\end{eqnarray*}

then for a fixed $c\in\Z^2\setminus \{0\}$, and $|n-m|>K$  with
$K\geq \frac{1}{\rho-\rho_+}\ln(\frac{\varepsilon}{\varepsilon_+})$,
$$\|\frac{\partial^2(P-R)}{\partial z_{n+tc}\partial \bar z_{m+tc}}-\lim_{t\to\infty}\frac{\partial^2(P-R)}
{\partial z_{n+tc}\partial \bar z_{m+tc}}\|\leq
\frac{\varepsilon}{|t|}e^{-|n-m|\rho} \leq
\frac{\varepsilon_+}{|t|}e^{-|n-m|\rho_+}.$$ That is  to say, $P-R$ satisfies $(A5)$ with $K_+,\varepsilon_+,\rho_+$ in place of $K,\varepsilon,\rho$.
 The proof of the remaining terms satisfying $(A5)$ is composed by the following two lemmas.

\begin{Lemma}\label{Ftoplitz}
$F$ satisfies $(A5)$ with $\varepsilon^{\frac 23}$ in place of $\varepsilon$.
\end{Lemma}

For the proof see \cite{GXY}.

\begin{Lemma}\label{toplitz}
Assume that $P$ satisfies $(A5)$, $F$ satisfies $(A5)$ with $\varepsilon^{\frac 23}$ in place of $\varepsilon$
and
$$\frac{\partial^2F}{\partial z_n\partial z_m}=0 (|n+m|>K),
\frac{\partial^2F}{\partial z_n\partial \bar z_m}=0
(|n-m|>K),\frac{\partial^2F}{\partial \bar z_n\partial {\bar z}_m}=0
(|n+m|>K),$$ then $\{P,F\}$ satisfies $(A6)$ with $\varepsilon_+$ in place of $\varepsilon$.
\end{Lemma}

For the proof see \cite{GXY}.\\

A KAM-step cycle is now completed.

\section{Iteration Lemma and Convergence}

\noindent

For any given $s,\varepsilon,r, \gamma$ and for all $\nu\ge 1$, we
define the following sequences \[
r_{\nu+1}=r(1-\sum_{i=2}^{\nu+2}2^{-i}),\]
\beq\varepsilon_{\nu+1}=c\gamma^{-5}(r_{\nu}-r_{\nu+1})^{-c}
K_{\nu}^{5(\tau+1)}\varepsilon_{\nu}^{\frac43},\eeq
\[\eta_{\nu+1}=\varepsilon_{\nu+1}^{\frac13}, L_{\nu+1}=L_{\nu}+\varepsilon_{\nu}\]
\[s_{\nu+1}=2^{-2}\eta_{\nu}s_{\nu}=2^{-2{(\nu+1)}}(\prod_{i=0}^{\nu}\varepsilon_i)^{\frac13}s_0,
\]
$$K_{\nu+1}^{1+3\varepsilon}\rho_{\nu+1}=1$$
$$K_{\nu+1}=3K_\nu=3^{\nu+1}K_0$$
\[\Delta_{\nu+1}=K^3_\nu
\]
where $c$ is a constant,$\gamma=\varepsilon_0^{\frac{1}{50}}\gg \varepsilon_0,$ and the parameters
$r_0,\varepsilon_0,s_0$ and $K_0$ are defined
 to be $r,\varepsilon,s$ and $K_0^2 e^{-K_0(r_0-r_1)}=\varepsilon_0^{\frac 13}$ respectively.

\subsection{Iteration lemma}

The preceding analysis can be summarized as follows.

\begin{Lemma}\label{Lem5.1}
Let $\varepsilon$ is small enough and $\nu\ge 0$. Suppose that

\noindent
(1). $N_\nu+{\cal B}_\nu+\bar{\cal B}_\nu$ is a normal form with parameters
$\xi$ satisfying \[ |\langle k,\omega_\nu\rangle|\ge \frac{\gamma}{K_\nu^\tau}, k\neq 0,\]
 \[|\langle k,\omega_\nu\rangle \pm\widetilde{\lambda}_j^\nu|\ge \frac{\gamma}{K_\nu^\tau},j\in{[n]} \]
 \[|\langle k,\omega_\nu\rangle \pm \widetilde{\lambda}_i^\nu\pm \widetilde{\lambda}_j^\nu|\ge \frac{\gamma}{K_\nu^\tau},i\in{[m]},j\in{[n]}\]
$$|det(\la k,\omega_\nu\ra I\pm{\cal A}_n^\nu\otimes I_2 \pm I_2\otimes {\cal A}_{n'}^\nu)|\geq \frac{\gamma}{K_\nu^\tau},k\neq 0,n,n'\in{\cal L}_2$$

 \beq\label{nonresonanceconditions}
 \eeq
\indent on a closed set $\Cal O_{\nu}$ of $\R^b$ for all $0<|k|\leq K_\nu$. Moreover,suppose that $\omega_\nu(\xi)$,$P^{011}_{ij\nu}(\xi)$,$A_{[n]}^{\nu}(\xi)$ are $C_W^4$ smooth and satisfy

$$\sup_{\xi\in \Cal O_\nu}\max_{d\leq 4}\|\partial_{\xi}^d(A_{[n]}^\nu-A_{[n]}^{\nu-1})\|\leq c\varepsilon_{\nu-1}$$
$$\sup_{\xi\in \Cal O_\nu}\max_{d\leq 4}|\partial_{\xi}^d(\Omega_{n'}^\nu-\Omega_{n'}^{\nu-1})|\leq \varepsilon_{\nu-1}$$
$$\sup_{\xi\in \Cal O_\nu}\max_{d\leq4}|\partial_\xi^d(\omega_\nu-\omega_{\nu-1})|\leq \varepsilon_{\nu-1}$$
and $$|P^{011}_{ij\nu}-P_{ij(\nu-1)}^{011}|_{\Cal O_\nu}\leq\varepsilon_{\nu-1} e^{-|i-j|\rho}$$
in the sense of Whitney.

\noindent (2).  $N_\nu+{\cal B}_\nu+\bar{\cal B}_\nu+P_\nu$ satisfies
$(A5)$ with $K_\nu,\varepsilon_\nu,\rho_\nu$ and
\[\|X_{P_\nu}\|_{D(r_\nu, s_\nu),\Cal O_{\nu}}\le
\varepsilon_\nu.\] Then there is a subset $\Cal
O_{\nu+1}\subset\Cal O_{\nu}$,
 \[\Cal O_{\nu+1}=\Cal
O_\nu\setminus(\Cal R_k^{\nu+1}),
\] $$\Cal R^{\nu+1}=\bigcup_{K_{\nu}<|k|\le K_{\nu+1},[n],[m],n,n'}(\Cal R_k^{\nu+1}\bigcup
\Cal R_{k[n]}^{\nu+1}\bigcup\Cal R_{k[n][m]}^{\nu+1}\bigcup{\Cal C^{\nu+1}_{knn'}(\gamma)}),
$$
where
 $$ \Cal R_k^{\nu+1}=\{\xi\in \Cal
O_{\nu}:|\langle k,\omega_{\nu+1}\rangle|< \frac{\gamma}{K_{\nu+1}^\tau}, k\neq 0\}$$
$$\Cal R_{k[n]}^{\nu+1}=\{\xi\in \Cal O_{\nu}:|\langle k,\omega_{\nu+1}\rangle \pm\widetilde{\lambda}_j^{\nu+1}|< \frac{\gamma}{K_{\nu+1}^\tau}\},$$
 $$\Cal R_{k[n][m]}^{\nu+1}=\{\xi\in \Cal O_{\nu}:\
|\langle k,\omega_{\nu+1}\rangle \pm \widetilde{\lambda}_i^{\nu+1}\pm \widetilde{\lambda}_j^{\nu+1}|< \frac{\gamma}{K_{\nu+1}^\tau},i\in [m],j\in[n]
\},$$
$$\Cal C^{\nu+1}_{knn'}=\{\xi\in \Cal O_{\nu}:\
|det(\la k,\omega_{\nu+1}\ra I\pm{\cal A}_n^{\nu+1}\otimes I_2 \pm I_2\otimes {\cal A}_{n'}^{\nu+1})|< \frac{\gamma}{K_{\nu+1}^\tau},k\neq 0,n,n'\in{\cal L}_2
\},$$

with $\omega_{\nu+1}=\omega_\nu+P_{0l00}^\nu$, and a symplectic
transformation of variables \beq \Phi_\nu:D_{\rho_\nu}(r_{\nu+1},s_{\nu +1})
\times\Cal O_{\nu}\to D_{\rho_\nu}(r_{\nu},s_{\nu}),\eeq such that on
$D_{\rho_{\nu+1}}(r_{\nu+1},s_{\nu +1})\times\Cal O_{\nu+1},$
$H_{\nu+1}=H_\nu\circ\Phi_\nu$ has the form
\begin{eqnarray*}
H_{\nu+1}&=&e_{\nu+1}+\la\omega_{\nu+1},I\ra+\sum_{[n]}
\la A_{[n]}^{\nu+1}(\xi)z_{[n]},\bar z_{[n]}\ra\nonumber\\
&&+\sum_{n'\in {\cal L}_2}[(\Omega_{n'}^{\nu+1}-\omega_{i'})z_{n'}\bar z_{n'}+(\Omega_{m'}^{\nu+1}-\omega_{j'})z_{m'}\bar z_{m'}]+{\cal B}_{\nu+1}+\bar{\cal B}_{\nu+1}+P_{\nu+1}.
 \end{eqnarray*}
 with
$$\sup_{\xi\in \Cal O_\nu}\max_{d\leq 4}\|\partial_{\xi}^d(A_{[n]}^{\nu+1}-A_{[n]}^{\nu})\|\leq c\varepsilon_{\nu}$$
$$\sup_{\xi\in \Cal O_\nu}\max_{d\leq 4}|\partial_{\xi}^d(\Omega_{n'}^{\nu+1}-\Omega_{n'}^{\nu})|\leq \varepsilon_{\nu}$$
$$\sup_{\xi\in \Cal O_\nu}\max_{d\leq4}|\partial_\xi^d(\omega_{\nu+1}-\omega_{\nu})|\leq \varepsilon_{\nu}$$
and $$|P^{011}_{ij(\nu+1)}-P_{ij\nu}^{011}|_{\Cal O_\nu}\leq\varepsilon_{\nu} e^{-|i-j|\rho}$$
in the sense of Whitney. And
$$\|X_{P_{\nu+1}}\|_{D(r_{\nu+1}, s_{\nu+1}),\Cal O_{{\nu+1}}}\le
\varepsilon_{\nu+1}.$$
\end{Lemma}

\subsection{Convergence}

Suppose that the assumptions of Theorem \ref{KAM} are satisfied to apply the iteration Lemma with $\nu=0$,recall that
$$\varepsilon_0=\varepsilon,r_0=r, s_0=s,L_0=L, N_0=N,{\cal B}_0={\cal B}, P_0=P,\gamma=\varepsilon^{\frac{1}{50}}, K_0^2 e^{-K_0(r_0-r_1)}=\varepsilon_0^{\frac 13}\quad
$$
$$
\Cal O_0= \left\{\xi\in \Cal O: \begin{array}{rcl} &&|\langle k,\omega\rangle|\ge \frac{\gamma}{K_0^\tau}, k\neq 0\\
 &&|\langle k,\omega\rangle \pm\widetilde{\lambda}_j|\ge \frac{\gamma}{K_0^\tau},j\in{[n]}\\
 &&|\langle k,\omega\rangle \pm \widetilde{\lambda}_i\pm \widetilde{\lambda}_j|\ge \frac{\gamma}{K_0^\tau},i\in{[m]},j\in{[n]}\\
&&|det(\la k,\omega\ra I\pm{\cal A}_n\otimes I_2 \pm I_2\otimes {\cal A}_{n'})|\geq \frac{\gamma}{K_0^\tau},n,n'\in{\cal L}_2
\end{array} \right\},$$
the assumptions of  the iteration lemma  are satisfied when
$\nu=0$ if $\varepsilon_0$ and $\gamma$ are sufficiently small.
Inductively, we obtain the following sequences:
\[
\Cal O_{\nu+1}\subset\Cal O_\nu,\]
\[\Psi^\nu=\Phi_0\circ\Phi_1\circ\cdots\circ\Phi_\nu:D_{\rho_\nu}(r_{\nu+1},s_{\nu+1})\times\Cal
O_\nu\to D_{\rho_0}(r_0,s_0),\nu\ge 0,
\]
\[H\circ\Psi^\nu=H_{\nu+1}=N_{\nu+1}+{\cal B}_{\nu+1}+\bar{\cal B}_{\nu+1}+P_{\nu+1}.\]
\indent Let $\tilde{\Cal O}=\cap_{\nu=0}^\infty \Cal O_\nu$. As in
\cite{P1,P2}, thanks to Lemma {\ref{Lem4.4}}, it concludes that
$N_\nu,\Psi^\nu,D\Psi^\nu,\omega_{\nu}$ converge uniformly on
$D_{\frac 12r}(\frac 12r,0)\times\tilde{\Cal O}$ with
 \begin{eqnarray*}
N_\infty+{\cal B}_\infty+\bar{\cal B}_\infty&=&e_\infty+\la\omega_\infty,I\ra+\sum_{[n]}
\la A_{[n]}^\infty(\xi)z_{[n]},\bar z_{[n]}\ra\nonumber\\
&&+\sum_{n'\in {\cal L}_2}[(\Omega_{n'}^{\infty}-\omega_{i'})z_{n'}\bar z_{n'}+(\Omega_{m'}^{\infty}-\omega_{j'})z_{m'}\bar z_{m'}]+{\cal B}_\infty+\bar{\cal B}_\infty.
 \end{eqnarray*}
\noindent Since
$$
\varepsilon_{\nu+1}=c\gamma^{-5}K_{\nu}^{5(\tau+1)}(r_\nu-r_{\nu-1})^{-c}\varepsilon_\nu^{\frac43}
, $$
  it follows that $\varepsilon_{\nu+1}\to 0$ provided that
$\varepsilon$ is sufficiently small. And we also have $\sum_{\nu=0}^{\infty}\varepsilon_{\nu}\leq2\varepsilon$.

Let $\phi_H^t$ be the flow of $X_H$. Since
$H\circ\Psi^\nu=H_{\nu+1}$, we have \beq\label{5.7}
\phi_H^t\circ\Psi^\nu=\Psi^\nu\circ\phi_{H_{\nu+1}}^t. \eeq The
uniform convergence of $\Psi^\nu,D\Psi^\nu,\omega_{\nu}$ and
$X_{H_{\nu}}$ implies that the limits can be taken on both sides
of (\ref{5.7}). Hence, on $D_{\frac 12r}(\frac 12r,0)\times\tilde{\Cal O}$ we
get \beq\label{5.8}
\phi_H^t\circ\Psi^\infty=\Psi^\infty\circ\phi_{H_{\infty}}^t\eeq
and
$$
\Psi^\infty:D_{\frac 12r}(\frac 12r,0)\times\tilde{\Cal O}\to
 D_\rho(r,s)  \times \Cal O.
 $$
  It follows from (\ref{5.8}) that
   $$
\phi_H^t(\Psi^\infty(\T^b\times
\{\xi\}))=\Psi^\infty\phi_{N_\infty}^t(\T^b\times\{\xi\})=\Psi^\infty(\T^b\times\{\xi\})
$$
\noindent for $\xi\in\tilde{\Cal O}$. This means that
$\Psi^\infty(\T^b\times\{\xi\})$ is an embedded torus which is
invariant for the original perturbed Hamiltonian system at $\xi\in
\tilde{\Cal O}$.  We remark here that the frequencies
$\omega_\infty(\xi)$ associated to
$\Psi^\infty(\T^b\times\{\xi\})$ are slightly different from
$\omega(\xi)$. The normal behavior of the invariant torus is
governed by normal frequencies $A_{[n]}^\infty,\Omega_{n'}^\infty$.\qed

\section{Measure Estimates}

This section is the essential part for this paper. For notational  convenience, let $\Cal O_{-1}=\Cal O$, $K_{-1}=0$.
Then at $\nu^{\rm th}$ step of  KAM iteration, we have to exclude
the following resonant set
$$\Cal R^{\nu+1}=\bigcup_{K_{\nu}<|k|\le K_{\nu+1},[n],[m],n,n'}(\Cal R_k^{\nu+1}\bigcup
\Cal R_{k[n]}^{\nu+1}\bigcup\Cal R_{k[n][m]}^{\nu+1}\bigcup{\Cal C^{\nu+1}_{knn'}(\gamma)}),
$$
where
 $$ \Cal R_k^{\nu+1}=\{\xi\in \Cal
O_{\nu}:|\langle k,\omega_{\nu+1}\rangle|< \frac{\gamma}{K_{\nu+1}^\tau}, k\neq 0\}$$
$$\Cal R_{k[n]}^{\nu+1}=\{\xi\in \Cal O_{\nu}:|\langle k,\omega_{\nu+1}\rangle \pm\widetilde{\lambda}_j^{\nu+1}|< \frac{\gamma}{K_{\nu+1}^\tau}\},$$
 $$\Cal R_{k[n][m]}^{\nu+1}=\{\xi\in \Cal O_{\nu}:\
|\langle k,\omega_{\nu+1}\rangle \pm \widetilde{\lambda}_i^{\nu+1}\pm \widetilde{\lambda}_j^{\nu+1}|< \frac{\gamma}{K_{\nu+1}^\tau},i\in [m],j\in[n]
\},$$
$$\Cal C^{\nu+1}_{knn'}=\{\xi\in \Cal O_{\nu}:\
|det(\la k,\omega_{\nu+1}\ra I\pm{\cal A}_n^{\nu+1}\otimes I_2 \pm I_2\otimes {\cal A}_{n'}^{\nu+1})|< \frac{\gamma}{K_{\nu+1}^\tau},k\neq 0,n,n'\in{\cal L}_2
\},$$

recall that $ \omega_{\nu+1}(\xi)=\omega(\xi)+\sum_{j=0}^\nu P_{0l00}(\xi)$ with $ |\sum_{j=0}^\nu P_{0l00}^{j}(\xi)|_{\Cal O_\nu}<\varepsilon $,and
$$\|A_{[n]}^{\nu +1}(\xi)-A_{[n]}(\xi)\|_{\Cal O_{\nu}}\leq \sum_{j=0}^\nu\|R_{[n][n]}^{011,j}\|\leq\varepsilon,$$
$$|\Omega_{n'}^{\nu +1}(\xi)-\Omega_{n'}(\xi)|_{\Cal O_{\nu}}\leq \sum_{j=0}^\nu |R_{n'n'}^{011,j}|\leq\varepsilon.$$

\noindent{\bf Remark.} From the section \ref{no4.4}, one has that
at $(\nu+1)^{\rm th}$ step, small divisor conditions are automatically
satisfied for $|k|\le K_{\nu}$. Hence, we only need to excise
the above resonant set $\Cal R^{\nu+1}$.

In the following, we only give the proof for the most complicated case $\{\xi\in \Cal
O_{\nu}:|\langle k,\omega_{\nu+1}\rangle + \widetilde{\lambda}_n^{\nu+1}- \widetilde{\lambda}_{n'}^{\nu+1}|< \frac{\gamma}{K_{\nu+1}^\tau},n,n'\in{\cal L}_1\}$ and $\{\xi\in \Cal O_{\nu}:|det(\la k,\omega_{\nu+1}\ra I+{\cal A}_n^{\nu+1}\otimes I_2 - I_2\otimes {\cal A}_{n'}^{\nu+1})|< \frac{\gamma}{K_{\nu+1}^\tau},n,n'\in{\cal L}_2\}$. When $n\in {\cal L}_1,n'\in {\cal L}_2$, there will be no small divisors. In other cases, the proof is similar, so we omit it. For simplicity, set $M^{\nu+1}=|\langle k,\omega_{\nu+1}\rangle + \widetilde{\lambda}_n^{\nu+1}- \widetilde{\lambda}_{n'}^{\nu+1}|$ and $Y^{\nu+1}=\la k,\omega_{\nu+1}\ra I+{\cal A}_n^{\nu+1}\otimes I_2 - I_2\otimes {\cal A}_{n'}^{\nu+1}$,$Y^{\nu}=\la k,\omega_{\nu}\ra I+{\cal A}_n^{\nu}\otimes I_2 - I_2\otimes {\cal A}_{n'}^{\nu}$,then for $|k|\leq K_\nu$
\begin{eqnarray*}
\|(Y^{\nu+1})^{-1}\|&=&\|(Y^\nu+(Y^{\nu+1}-Y^{\nu}))^{-1}\|\\
&=&\|(I+(Y^{\nu})^{-1}(Y^{\nu+1}-Y^{\nu}))^{-1}(Y^{\nu})^{-1}\|\\
&\leq&2\|(Y^{\nu})^{-1}\|\leq2\frac{K^\tau_{\nu}}{\gamma}\leq\frac{K^\tau_{\nu+1}}{\gamma}.
\end{eqnarray*}

\begin{Lemma}
For any given $n, n'\in \Bbb Z_1^2$ with $|n-n'|\leq K_{\nu+1}$, either $|\langle k,\omega_{\nu+1}\rangle + \widetilde{\lambda}_n^{\nu+1}- \widetilde{\lambda}_{n'}^{\nu+1}|>1$ or there are $n_0, n'_0, c\in\Z^2$ with $|n_0|, |n'_0|, |c|\le  3K_{\nu+1}^2$ and $t\in  \Z$, such that $n=n_0+tc$, $n'=n'_0+tc$.
\end{Lemma}

\proof Since $|n-n'|\leq K_{\nu+1}$,  with an elementary calculation
\begin{eqnarray*}|n|^2-|n'|^2=|n-n'|^2+2\langle
n-n',n'\rangle
\end{eqnarray*}
If $|\langle n-n',n'\rangle|>K_{\nu+1}^2$,
we have $|\langle k,\omega_{\nu+1}\rangle + \widetilde{\lambda}_n^{\nu+1}- \widetilde{\lambda}_{n'}^{\nu+1}|>1$, there will be no small divisor.

In the case that $|\langle n-n',n'\rangle|\le K_{\nu+1}^2 $, clearly $n-n'=0$ is trivial. Assume $n-n'\neq 0$, without loss of generality, we assume that the first component $(n-n')_1$ of $n-n'$ is not zero. Let

$$c=(-(n-n')_2,(n-n')_1)$$

Then
$$c\perp (n-n')$$
and $c\in\Z^2\setminus\{0\}$ with $|c|\le |n-n'|\le K_{\nu+1}$. Clearly, $ c, n-n'$ are linearly independent, hence there exist $x_1, x_2\in \R$ such that
$$n'=x_1c+x_2(n-n').$$
 Set (here $[\cdot]$ denotes the integer part of $\cdot$)
 $$t=[x_1]$$
 then $t\in\Z$ and $|n'-tc|\le 2K_{\nu+1}^2$.
Take $n'_0= n'-tc\in\Z^2$ and $n_0=n'_0+n-n'\in\Z^2$. We have $|n'_0|\le 2K_{\nu+1}^2$ and
\[ |n_0|\le |n'_0|+|n-n'|\le 3K_{\nu+1}^2.\]
\qed

\begin{Lemma}\[\cup_{n,n'\in {\cal L}_1} \Cal R^{\nu+1}_{k[n][n']}\subset
\cup_{{n_0,n'_0,c}\in \Bbb Z^2, t\in\Z} \Cal R_{k, n_0+tc,
n'_0+tc}^{\nu+1}\] where $|n_0|, |n'_0|, |c|\le  3K_{\nu+1}^2$.
\end{Lemma}

\proof If $|\langle n-n',n'\rangle|>
 K_{\nu+1}^2$, $\Cal R^{\nu+1}_{k[n][n']}=\emptyset.$  If $|\langle n-n',n' \rangle|\le K_{\nu+1}^2$,
there exist $n_0, n'_0, c\in \Bbb Z^2, t\in\Z$ with $|n_0|, |n'_0|, |c|\le  3K_{\nu+1}^2$ such
that $n=n_0+tc$, $n'=n'_0+tc$. Hence \[\cup_{n,n'\in {\cal L}_1} \Cal R^{\nu+1}_{k[n][n']}\subset
\cup_{{n_0,n'_0,c}\in \Bbb Z^2, t\in\Z} \Cal R_{k, n_0+tc,
n'_0+tc}^{\nu+1}\] where $|n_0|, |n'_0|, |c|\le  3K_{\nu+1}^2$.\qed

\begin{Lemma} For fixed $k,n_0, n'_0, c$, one has
 $$
{\rm meas}(\cup_{t\in\Z} \Cal R_{k, n_0+tc, n'_0+tc}^{\nu+1})<
c\frac{\gamma}{K_{\nu+1}^{\tau\over{2} }}.
$$
\end{Lemma}

\proof
Due to T\"oplitz-Lipschitz property of
 $N_{\nu}+{\cal B}_{\nu}+\bar{\cal B}_{\nu}+P_{\nu}$,
then
 $$|M^{\nu+1}(t)-\lim_{t\to\infty}M^{\nu+1}(t) |< \frac {\varepsilon_0}{|t|}.$$

We define resonant set \beq\label{resonant3} \Cal
R_{kn_0n'_0c\infty}^{\nu+1}=\{\xi\in \Cal
O_{\nu}:|\lim_{t\to\infty}M^{\nu+1}(t))|<\frac
{\gamma}{K_{\nu+1}^{\tau\over{2}}}\}\nonumber\eeq

For fixed $k,n_0, n'_0, c$,
$${\rm meas}(\Cal R_{kn_0n'_0c\infty}^{\nu+1})<
\frac{\gamma}{K_{\nu+1}^{\tau\over{2}}}.
$$
Then for   $\xi\in\Cal O_{\nu}\backslash \Cal
R_{kn_0n'_0c\infty}^{\nu+1}$, we have
$$|\lim_{t\to\infty}M^{\nu+1}(t))|\ge
\frac{\gamma}{K_{\nu+1}^{\tau\over{2}}}.$$

Case 1: When $|t|>K_{\nu+1}^{\tau\over{2}}$,  for   $\xi\in\Cal O_{\nu}\backslash \Cal
R_{kn_0n'_0c\infty}^{\nu+1}$, we have

\begin{eqnarray*}
&&|M^{\nu+1}(t)|\\
&\geq&|\lim_{t\to\infty}M^{\nu+1}(t)|-\frac {\varepsilon_0}{|t|}\\
&\geq& \frac
{\gamma}{K_{\nu+1}^{\tau\over{2}}}-\frac
{\varepsilon_0}{K_{\nu+1}^{\tau\over{2}}}\\
&\geq& \frac
{\gamma}{2K_{\nu+1}^{\tau\over{2}}}.\end{eqnarray*}

Case 2: When $|t_{1}|\le K_{\nu+1}^{\tau\over{2}}$, we define resonant set \beq \Cal
R_{kn_0n'_0ct}^{\nu+1}=\{\xi\in \Cal
O_{\nu}:|M^{\nu+1}(t)|<\frac
{\gamma}{K_{\nu+1}^{\tau}}\}\nonumber\eeq

For fixed $k,n_0, n'_0, c,t$,
$${\rm meas}(\Cal R_{kn_0n'_0ct}^{\nu+1})<
\frac{\gamma}{K_{\nu+1}^{\tau}},
$$
then
$${\rm meas}\{\cup_{|t|\le K_{\nu+1}^{\tau\over{2}}}\Cal R_{kn_0n'_0ct}^{\nu+1}\}<K_{\nu+1}^{{\tau\over{2}}}\frac
{\gamma}{K_{\nu+1}^{\tau}}\le \frac {\gamma}{K_{\nu+1}^{\tau\over{2}}}.$$

As a consequence,
$$
{\rm meas}(\cup_{t\in\Z} \Cal R_{k, n_0+tc, n'_0+tc}^{\nu+1})<
c\frac{\gamma}{K_{\nu+1}^{\tau\over{2} }}.
$$\qed

For $K_\nu<|k|\leq K_{\nu+1}$,we consider $n,n'\in {\cal L}_2$ as an example,the other cases can be proved analogously.Assume that $(n,m)$ and $(n',m')$ are resonant pairs in ${\cal L}_2$,then

\begin{Lemma}
For any given $n, n'\in \Bbb Z_1^2$ with $|n-n'|\leq K_{\nu+1}$, either $|det(\la k,\omega_{\nu+1}\ra I+{\cal A}_n^{\nu+1}\otimes I_2 - I_2\otimes {\cal A}_{n'}^{\nu+1})|>1$ or there are $n_0, n'_0, c\in\Z^2$ with $|n_0|, |n'_0|, |c|\le  3K_{\nu+1}^2$ and $t\in  \Z$, such that $n=n_0+tc$, $n'=n'_0+tc$.
\end{Lemma}

\begin{Lemma}
\[\cup_{n,n'\in \Bbb Z_1^2} \Cal C_{knn'}^{\nu+1}\subset
\cup_{{n_0,n'_0,c}\in \Bbb Z^2, t\in\Z} \Cal C_{k, n_0+tc,
n'_0+tc}^{\nu+1}\] where $|n_0|, |n'_0|, |c|\le  3K_{\nu+1}^2$.
\end{Lemma}

\begin{Lemma} For fixed $k,n_0, n'_0, c$, one has
 $$
{\rm meas}(\cup_{t\in\Z} \Cal C_{k, n_0+tc, n'_0+tc}^{\nu+1})<
c\frac{\gamma^{\frac14}}{K_{\nu+1}^{\frac{\tau}{20}}}.
$$
\end{Lemma}
\proof
Due to the analysis above and T\"{o}plitz-Lipschitz property of $N+{\cal B}+\bar{\cal B}+P$,the coefficient matrix $Y^{\nu+1}(t)$ has a limit as $t\rightarrow\infty$,
$$\|Y^{\nu+1}(t)-\lim_{t\rightarrow\infty}Y^{\nu+1}(t)\|\leq \frac{\varepsilon_0}{t}.$$
We define resonant set
$$\Cal C_{kn_0n'_0c\infty}^{\nu+1}=\left\{\xi\in \Cal O_\nu:|det\lim_{t\rightarrow\infty}Y^{\nu+1}(t)|< \frac{\gamma}{K_{\nu+1}^{\frac{\tau}{5}}}\right\}.$$

Then for   $\xi\in\Cal O_{\nu}\backslash \Cal
C_{kn_0n'_0c\infty}^{\nu+1}$, we have
$$\|(\lim_{t\rightarrow\infty}Y^{\nu+1}(t))^{-1}\|\leq
\frac{K_{\nu+1}^{\frac{\tau}{5}}}{\gamma}.$$

Since
$$\|Y^{\nu+1}(t)-\lim_{t\rightarrow\infty}Y^{\nu+1}(t)\|\leq \frac{\varepsilon_0}{t},$$
for $|t|>K_{\nu+1}^{\frac{\tau}{5}}$,we have
$$\|(Y^{\nu+1}(t))^{-1}\|\leq 2\frac{K_{\nu+1}^{\frac{\tau}{5}}}{\gamma}\leq\frac{K_{\nu+1}^{\tau}}{\gamma}.$$

For $|t|\leq K_{\nu+1}^{\frac{\tau}{5}}$, we define resonant set \beq \Cal
C_{kn_0n'_0ct}^{\nu+1}=\{\xi\in \Cal
O_{\nu}:|detY^{\nu+1}(t)|< \frac{\gamma}{K_{\nu+1}^{\tau}}\}.\nonumber\eeq

In addition
$$\inf_{\xi\in \Cal O}\max_{0<d\leq 4}|\partial^d_\xi(detY^{\nu+1}(t))|\geq\frac12|k|\geq\frac12K.$$

For fixed $k,n_0, n'_0, c,t$,
$${\rm meas}(\Cal C_{kn_0n'_0ct}^{\nu+1})<
(\frac{\gamma}{K_{\nu+1}^{\tau}})^{\frac14},
$$
then
$${\rm meas}\{\cup_{|t|\le K_{\nu+1}^{\tau\over{5}}}\Cal C_{kn_0n'_0ct}^{\nu+1}\}<K_{\nu+1}^{{\tau\over{5}}}(\frac
{\gamma}{K_{\nu+1}^{\tau}})^{\frac14}\le \frac {\gamma^{\frac14}}{K_{\nu+1}^{{\tau\over{20}}}}.$$

As a consequence,
$$
{\rm meas}(\cup_{t\in\Z} \Cal C_{k, n_0+tc, n'_0+tc}^{\nu+1})<
c\frac {\gamma^{\frac14}}{K_{\nu+1}^{{\tau\over{20}}}}.
$$\qed

\begin{Lemma}\label{Lem6.1}
$${\rm meas}(\bigcup_{K_{\nu}<|k|\le K_{\nu+1}} R_k^{\nu+1})\leq cK_{\nu+1}^{b}\frac {\gamma}{K_{\nu+1}^{\tau}}=c\frac {\gamma}{K_{\nu+1}^{\tau-b}}$$

$${\rm meas}(\bigcup_{K_{\nu}<|k|\le K_{\nu+1},[n]} R_{k[n]}^\nu)\leq cK_{\nu+1}^{2+b}\frac {\gamma}{K_{\nu+1}^{\tau}}=c\frac {\gamma}{K_{\nu+1}^{\tau-2-b}}$$

$${\rm meas}(\bigcup_{K_{\nu}<|k|\le K_{\nu+1},[n],[ m]} R_{k[n][m]}^{\nu+1})\leq  c\frac {\gamma}{K_{\nu+1}^{{\tau\over{2}}-12-b}}$$

$${\rm meas}(\bigcup_{K_{\nu}<|k|\le K_{\nu+1},n,n'} \Cal C_{knn'}^\nu)\leq c\frac {\gamma^{\frac14}}{K_{\nu+1}^{{\frac{\tau}{20}}-12-b}}$$
\end{Lemma}

\begin{Lemma}\label{Lem6.2} Let $\tau>20(12+b+1)$, then
the total measure need to exclude along the KAM iteration is
\begin{eqnarray*}
&&{\rm meas}(\bigcup_{\nu\ge 0}\Cal R^{\nu+1})\\
 &=&{\rm
meas}[\bigcup_{\nu\ge 0}\bigcup_{K_{\nu}<|k|\le
K_{\nu+1},[n],[m],n,n'}(\Cal R_k^{\nu+1}\bigcup
\Cal R_{k[n]}^{\nu+1}\bigcup\Cal R_{k[n][m]}^{\nu+1}\bigcup{\Cal C^{\nu+1}_{knn'}(\gamma)})]\\
&\le&c\sum_{\nu\ge 0}\frac{\gamma^{\frac14}}{K_{\nu+1}}\le c\gamma^{\frac14}.
\end{eqnarray*}
\end{Lemma}

\section{Appendix}

\sss
\begin{Lemma}\label{Lem7.1}
$$\|FG\|_{ D_\rho(r,s),\Cal O  }\le \|F\|_{ D_\rho(r,s),\Cal O  }\|G\|_{ D_\rho(r,s),\Cal O  }.$$\end{Lemma}

\proof Since
$(FG)_{kl\alpha\beta}=\sum_{k',l',\alpha',\beta'}F_{k-k',l-l',\alpha-\alpha',\beta-\beta'}
G_{k'l'\alpha'\beta'}$, we have
\begin{eqnarray*}
\|FG\|_{ D_\rho(r,s),\Cal O
}&=&\sup_{D_\rho(r,s)}\sum_{k,l,\alpha,\beta}|(FG)_{kl\alpha\beta}|_{\Cal
O}|I^{l}||z^{\alpha}| |\bar
z^{\beta}|e^{|k||{\rm Im}\theta|}\\
&\le&\sup_{D_\rho(r,s)}
\sum_{k,l,\alpha,\beta}\sum_{k',l',\alpha',\beta'}|F_{k-k',l-l',\alpha-\alpha',\beta-\beta'
}G_{k'l'\alpha'\beta'}|_{\Cal O}|I^{l}||z^{\alpha}| |\bar
z^{\beta}|e^{|k||{\rm Im}\theta|}\\
&\le&\|F\|_{ D_\rho(r,s),\Cal O  }\|G\|_{ D_\rho(r,s),\Cal O  }
\end{eqnarray*}
and the proof is finished.\qed

\begin{Lemma}\label{Lem7.2}
(Generalized Cauchy inequalities)
$$
\|F_{\theta}\|_{D_\rho(r-\sigma,s),\Cal O}\le
\frac{c}{\sigma}\|F\|_{ D_\rho(r,s),\Cal O },$$
$$
\|F_{I}\|_{D_\rho(r,\frac 12 s),\Cal O}\le \frac {c}{s^2}\|F\|_{
D_\rho(r,s),\Cal O },$$ and
$$
\|F_{z}\|_{D_{\rho}(r,\frac 12 s),\Cal O}\le
\frac{c}{s}\|F\|_{ D_\rho(r,s),\Cal O },$$
$$
\|F_{\bar z}\|_{D_{\rho}(r,\frac 12 s),\Cal O}\le
\frac{c}{s}\|F\|_{ D_\rho(r,s),\Cal O }.$$
\end{Lemma}
\proof We only prove the third inequality, the others can be proved
similarly. Let $w\neq 0$, then $f(t)=F(z+tw)$ is an analytic map
from the complex disc $|t|<\frac{s}{\|w\|_\rho}$ in $\C$ into
$D_\rho(r,s)$. Hence $$\|f'(0)\|_{D_\rho(r,\frac 12 s),\Cal O}=\|
F_zw\|_{D_\rho(r,\frac 12 s),\Cal O}\leq \frac{c}{s}\|F\|_{
D_\rho(r,s),\Cal O }\cdot \|w\|_\rho,$$ by the usual Cauchy
inequality. Since $w\neq 0$, so
$$\frac{\|F_zw\|_{D_\rho(r,\frac 12 s),\Cal
O}}{\|w\|_\rho}\leq \frac{c}{s}\|F\|_{ D_\rho(r,s),\Cal O },$$ thus
$$\|F_{z}\|_{{D_\rho(r,\frac 12 s),\Cal O}}=\sup_{w\neq 0}\frac{\|F_zw\|_{D_\rho(r,\frac 12 s),\Cal O}}{\|w\|_\rho}\leq
\frac{c}{s}\|F\|_{ D_\rho(r,s),\Cal O }.$$ \qed

 \indent Let
$\{\cdot,\cdot\}$ denote the Poisson bracket of smooth functions,
i.e.,
\[ \{F,G\}=\la\frac{\partial F}{\partial
I}, \frac{\partial G}{\partial \theta}\ra-\la \frac{\partial F}{
\partial \theta},\frac{\partial G}{\partial I}\ra+{\rm i}(
\langle\frac{\partial F}{\partial z},\frac{\partial G} {\partial
{\bar z}}\rangle-\langle \frac{\partial F}{\partial {\bar
z}},\frac{\partial G} {\partial { z}}\rangle),\] then we have the
following lemma:

\begin{Lemma}\label{Lem7.3} If
$$ \|X_F\|_{ D_\rho(r,s),\Cal O  }< \varepsilon',\ \|X_G\|_{ D_\rho(r,s),\Cal O  }< \varepsilon'', $$
then
$$
\|X_{\{F,G\}}\|_{D_\rho(r-\sigma,\eta s),\Cal
O}<c\sigma^{-1}\eta^{-2}\varepsilon'\varepsilon'',\ \eta\ll 1.$$ In
particular,  if $\eta\sim\varepsilon^{\frac 13}$, $\varepsilon',
\varepsilon''\sim \varepsilon$, we have
$\|X_{\{F,G\}}\|_{D_\rho(r-\sigma,\eta s),\Cal O}\sim
\varepsilon^{\frac {4}{3}}$.
\end{Lemma}
\proof By Lemma \ref{Lem7.1} and Lemma \ref{Lem7.2},
\begin{eqnarray*}
\|\frac{\partial^2F}{\partial I\partial I}\frac{\partial G}{\partial
\theta}\|_{D_\rho(r-\sigma,\frac{1}{2}s)}&<&
c\sigma^{-1}s^{-2}\|\frac{\partial F}{\partial I}\|_{D_\rho(r,s)}\cdot\|\frac{\partial G}{\partial \theta}\|_{D_\rho(r,s)},\\
\|\frac{\partial^2F}{\partial I\partial \theta}\frac{\partial
G}{\partial \theta}\|_{D_\rho(r-\sigma,\frac{1}{2}s)}&<&
c\sigma^{-1}\|\frac{\partial F}{\partial
I}\|_{D_\rho(r,s)}\cdot\|\frac{\partial G}{\partial
\theta}\|_{D_\rho(r,s)},\\
\|\frac{\partial^2F}{\partial I\partial z}\frac{\partial G}{\partial
\theta}\|_{D_\rho(r-\sigma,\frac{1}{2}s)}&<&
c\sigma^{-1}s^{-1}\|\frac{\partial F}{\partial
I}\|_{D_\rho(r,s)}\cdot\|\frac{\partial G}{\partial
\theta}\|_{D_\rho(r,s)},\\
\|\frac{\partial^2F}{\partial I\partial \bar z}\frac{\partial
G}{\partial \theta}\|_{D_\rho(r-\sigma,\frac{1}{2}s)}&<&
c\sigma^{-1}s^{-1}\|\frac{\partial F}{\partial
I}\|_{D_\rho(r,s)}\cdot\|\frac{\partial G}{\partial
\theta}\|_{D_\rho(r,s)},\\
\|\frac{\partial^2F}{\partial z\partial I}\frac{\partial G}{\partial
\bar z}\|_{D_\rho(r-\sigma,\frac{1}{2}s)}&<&
c\sigma^{-1}s^{-2}\|\frac{\partial F}{\partial
z}\|_{D_\rho(r,s)}\cdot\|\frac{\partial G}{\partial \bar
z}\|_{D_\rho(r,s)},\\
\|\frac{\partial^2F}{\partial z\partial \theta}\frac{\partial
G}{\partial \bar z}\|_{D_\rho(r-\sigma,\frac{1}{2}s)}&<&
c\sigma^{-1}\|\frac{\partial F}{\partial
z}\|_{D_\rho(r,s)}\cdot\|\frac{\partial G}{\partial \bar
z}\|_{D_\rho(r,s)},\\
\|\frac{\partial^2F}{\partial z\partial z}\frac{\partial G}{\partial
\bar z}\|_{D_\rho(r-\sigma,\frac{1}{2}s)}&<&
c\sigma^{-1}s^{-1}\|\frac{\partial F}{\partial
z}\|_{D_\rho(r,s)}\cdot\|\frac{\partial G}{\partial \bar
z}\|_{D_\rho(r,s)},\\
\|\frac{\partial^2F}{\partial z\partial \bar z}\frac{\partial
G}{\partial \bar z}\|_{D_\rho(r-\sigma,\frac{1}{2}s)}&<&
c\sigma^{-1}s^{-1}\|\frac{\partial F}{\partial
z}\|_{D_\rho(r,s)}\cdot\|\frac{\partial G}{\partial \bar
z}\|_{D_\rho(r,s)}.
\end{eqnarray*}
The other cases can be obtained analogously, hence
$$
\|X_{\{F,G\}}\|_{D_\rho(r-\sigma,\eta s),\Cal
O}<c\sigma^{-1}\eta^{-2}\varepsilon'\varepsilon''.
$$\qed

\

\

  \end{document}